# Presto! Digitization:  Part I

## From NKS Number Theory to "XORbitant" Semantics, by way of Cayley-Dickson Process and Zero-Divisor-based "Representations"

Robert de Marrais
rdemarrais@alum.mit.edu

## ABSTRACT

The objects of the great Nonlinear Revolutions – Catastrophes and Chaos in the 1960s-70s (henceforth, CT); and, "small-world" and "scale-free" Network Theory (NT), emerging from studies of synchronization since the Millennium – will be spliced together by a New Kind of Number Theory, focused on digitizations (i.e., binary strings).  NT nodes then become "feature-rich" representations (nodules in a "rhizosphere") of CT contents.

The "Box-Kite" formalism of zero divisors (ZD's) – first showing in the 16-D Sedenions, then in all higher $2^N$-ions derived from Imaginaries by Cayley-Dickson Process (CDP) – can model such "enriched" nodules.  Its (bit-string XOR-ing vs. matrix-multiplication-based) operations unfold into "representations" of the objects traditionally linked with "partitions of Nullity":  *Singularities*.  The route from this "local" level of CT to fractals and Chaos, via CDP extensions to $2^N$-ions for growing N, will involve us in graphics of higher-dimensional "carry-bit overflow," as manifest in the mandala-like patterns showing up in "emanation tables" (the rough equivalent, for ZD's, of group theorists' Cayley Tables).  More concretely, this route will model semantic propagation – "meaning" itself.

I'll lead into this with a quote about "Hjelmslev's Net" (which I'll claim is the CDP *manqué*) from a famous postmodern text deeply engaged with the philosophy of future mathematics, Deleuze and Guattari's A Thousand Plateaus (where "rhizosphere" imagery arose).  From here, with strong assists from the CT-based structuralism of Jean Petitot, Algirdas Greimas's "Semiotic Square" will show us how to explicitly link CT thinking to the foundations of semantics via ZD "representations" (Part III), while the infinite-dimensional ZD meta-fractal or "Sky" which Box-Kites fly "beneath" or "in" – first appearing in the 32-D Pathions and incorporating the higher $2^N$-ions, and the theme of Part II – will provide sufficient *lebensraum* for Lévi-Strauss's "Canonical Formula" of mythwork to unfurl in, thereby doing for semiotics what S-matrix methods promised for particle physics. (These results serve to extend my NKS 2004 paper, available online at the Wolfram Science website.)

## 0. Deep Background:  Size vs. Place in the History of Number

The ingenious method of expressing every possible number using a set of ten symbols (each symbol having a place value and an absolute value) emerged in India. The idea seems so simple nowadays that its significance and profound importance is no longer appreciated. Its simplicity lies in the way it facilitated calculation and placed arithmetic foremost amongst useful inventions. The importance of this invention is more readily appreciated when one considers that it was beyond the two greatest men of Antiquity, Archimedes and Apollonius.

*Pierre-Simon, Marquis de Laplace*[1]

From the vantage of traditional Number Theory, based on Quantity and Measure, and a notion of "Natural" synonymous with indices of Counting, what could be more *exorbitant* – etymologically, more "off the beaten path" – than basing Number, instead, on the "XOR bit-dance" implicit in strings of binary toggle-switches?  The notion of Zero was unknown to the Greeks, among whom Number Theory first flourished.  It reached the West from China, via India, then Arabic North Africa, at the dawn of the Italian Renaissance, where the powerful utility of its "place-holder" thinking (especially when compared to the clunkiness of Roman-numeral manipulations) impressed the new breed of double-entry bookkeeping-savvy accountants – in particular, Leonardo of Pisa, son of an Algerian trading post's director, best known by his nickname, Fibonacci.  His <u>Book of the Abacus</u> of 1202 introduced Hindu-Arabic computation, applying it to many concrete situations (one of which, involving penned rabbits, led to the series since named for him).

In his later, deeper but less famous <u>Book of Squares</u>, Fibonacci proved himself the greatest Number Theorist between Diophantus of mid-3[rd] Century Alexandria and Pierre de Fermat of 17[th] Century France.  In particular, his extension of classical results concerning Pythagorean triplets (integer sides of right triangles) led the way toward 3 great modern insights.  First, Hilbert's generalization of what a metric is, in a geometry exceeding (in a manner even more extreme than Descartes' algebraization anticipated) the 3-D limits of Greek thinking:  its summing an infinite number of squares (mediated by coordinates-based right triangles), to yield one more (the square of the "distance"), is the bedrock of classical Quantum Mechanics (QM).  Second, by positing the uniqueness of *squaring* (as opposed to cubing or other empowering) of integer trios, in the manner of the formula ascribed to (but much older than) Pythagoras, Fermat's Last Theorem of 1630 presented mathematicians with their most profound (yet simply statable) conundrum – one which would trigger momentous developments like Kummer's theory of ideals, and serve (up through the work of Andrew Wiles) as a guide to pure research for some four centuries.  But third, the number theoretics epitomized by Fibonacci's study of triplets finally ran aground, in Hurwitz's Proof a century ago, when attempts to assimilate the second, "calendrical" sense of Number as *Cycle*, to the "linear" notion based on *Counting*, came unraveled – when generalizations, that is, of the *angular* sense of exponentiation, typifying Imaginary numbers in 2-D, broke down at last with their third dimension-doubling, by the Cayley-Dickson Process (CDP), to the 16-D Sedenions . . . where the notion  of "distance" itself begins to totter, and Zero-Divisors appear "out of nowhere."

In the domain of pure Number, the failure of metrics and emergence of zero-divisors (ZD's) are joined together, but have their own flavors. There is a progressive breakdown of other properties at each dimension-doubling: the multiplicity of "sheets" associated by Riemann with equations over fields of Imaginary Number (a "parking garage" with infinite levels for functions as fundamental as the logarithmic) is the immediate side-effect of angular exponents, since there will always be an infinite number of angles with the same value but different "winding number" (just add $\pm 2\pi k$, k any integer, to the exponent's "angle"). In the 4-D Quaternions, we lose commutativity; with the 8-D Octonions, we still have a "field," but it's no longer associative. With each of these losses of structure, however, there are clear-cut gains in sensitivity: when equations are permitted to have Complex solutions, the fundamental theorem of algebra – linking precisely N (not necessarily distinct) solutions or "Zero's" with any $N^{th}$ order polynomial – becomes almost trivial. With Quaternions, the vector calculus of "dot" and "cross" products Gibbs extracted from them makes the intricacies of physical mechanics simple enough for college freshmen – with the "full Quaternion," meanwhile, serving to underwrite QM's angular momentum operators and Uncertainty relations, as well as the first nonlinear (and "non-Abelian," hence at least Quaternionic) gauge theories.

The 16-D loss of a metric has its benefits as well: I'll argue that, in fact, it provides the "natural" context within which we can think about the basic intuition of "small-world" networks, wherein the "strength of weak ties"[2] is paramount: for here, the key is gauging relations not of kinship, but acquaintance –"six degrees of Kevin Bacon," not closed systems of blood relationship. Here, if "~" means "is acquainted with," we must expect that if A ~ B and B ~ C, then C ~ A in general does *not* obtain, which is to say the basis of a job-hunter's "networking" does *not* conform to a standard metric. The metrification of $2^N$-ions, N < 4, is the straightforward result of their all deploying generalized "Pythagorean squares" patterns: writing an arbitrary complex number as $z = a + bi$, and another as $z^* = c + di$, we take it for granted that the product of their absolute values is the absolute value of their products (i.e., it "has a norm"); alternatively, writing things out in terms of components, we get an extended sort of Pythagorean relationship between squares: the product of *one* sum of N squares by *another* is a sum of N squares as well, with N = 2, 4, or 8, for Complex, Quaternionic, or Octonionic Numbers respectively. In the Complex case, this just means $(ac - bd)^2 + (ad + bc)^2 = (a^2 + b^2)(c^2 + d^2).$

What Hurwitz showed in 1898 is that such patterns no longer work for N > 8. If, using Rajwade's shorthand,[3] we rewrite the rules just given to have form (N, N, N), where the first two N's stand for the number of squares in the multiplied terms, and the final N for those in their product, then the Sedenions conform "at worst" to the rule (16, 16, 32): since 16 = 8 + 8, apply the Octonions' "8-squares" rule 4 times. ("At worst" means the final index in parentheses, at the time Rajwade was writing, *may* be < 32, although it was known at the time to be *at least* 29.) Indefinite redoubling yields a "worst case" rule (and harbinger of infinite-dimensional ZD "meta-fractals" we'll call "skies") for $2^N$-ions, N = $(3 + k) \geq 4$: $(2^N, 2^N, 2^{N+k}).$ This is a side-effect of "Bott periodicity," an 8-D modularity found in countless arenas (e.g., Clifford algebras and cobordism theories) since Raoul Bott first displayed it in more limited topological contexts. The keyword is *modularity*: for ZD's arise simply in ring theory as soon as we consider *finite* rings, over integer domains made from *composites*, not primes. Consider the integers modulo 6: 1 + 5 = 2 + 4 = 3 + 3 = 0, but 2 * 3 = 0 as well, making 2 and 3 "mutual ZD's" in this context.

It was hardly obvious, though, that ZD's were an inevitable side-effect of the modularity of angular exponents. Gauss himself had assumed the Imaginary realm was the "shadow of shadows" of an infinite hierarchy of higher entities. (The "composite" character – not of $2^N$-ions for some high N, but of Complex Numbers themselves – was a shock. Kummer's theory of ideals arose in response to the collapse of his "proof" of Fermat's Last Theorem, which foundered on the non-unique nature of decomposing Complex "number rings" into prime factors.) When plausible extensions of the Imaginary were found to occur only in $2^N$ dimensions, and – as Hamilton first showed with the failure of the commutative law among Quaternions – further displayed surprising losses of just those behaviors thought indispensable to any sane notion of Number, the suspicion that things were not so simple as Gauss had supposed started growing, and eventually was confirmed by Hurwitz – after which, CDP-generated hypercomplex numbers were discarded as "pathological." (A situation highly reminiscent of that which Benoit Mandelbrot found in analysis, whose "monsters" he turned into the "pets" we now call "fractals." The connection between these, meanwhile, and the XOR-based "carry-bit overflow" that characterizes progressive derivations of $2^N$-ions by CDP, will provide the bridge we need between ZD's and "scale-free" – i.e., "fractal power-law" – networks.)

And so the situation remained for exactly a century, until Guillermo Moreno's 1998 paper[4] indicated ZD's conformed to some remarkably simple patterns (they span a subspace of the Sedenions homomorphic to $G_2$, the Octonions' automorphism group – a result extensible, in turn, to all higher Number forms derived by CDP, since $G_2$ was known already to be the general $2^N$-ions' "derivation algebra" from $2^{N-1}$-ions). The discovery, 2 years later, that ZD's obeyed rules of engagement most cleary displayed in the workings of XOR relations between bit-strings – that is, by an NKS, not classical, sense of what Number is – marked the first results of my own research. (In a twisted echo of N-squares norms, we'll find our "master key" in this: the index of the product of two $2^N$-ions is the XOR of their indices!) But the hundred years between Hurwitz and the new interest in the "metric-free" ways of ZD's didn't see the problem *dismissed* so much as *displaced*. As my argument hinges on this displacement, let's give it a special name.

Lévi-Strauss has famously compared the structures of Baroque and Classical music to those his analyses found in North and South American mythologies, explaining this correspondence as a kind of "reincarnation."[5] As the advent of printing presses and print-specific modes of rational expression (viz., novels) replaced mythic narrative, the *structures* of the latter were reborn in new forms, in the "well-tempered'" modes of music we associate with names like J. S. Bach, Haydn, and Mozart. With the advent of communications technology and electronic modes of broadcasting and conversing, a similar (but oppositely directed) "reincarnation" can be indicated: the program killed by Hurwitz's proof in mathematics is reborn, ensouled in the new body of *structural linguistics*.

The genealogical derivation of Number Forms via CDP – which, starting with the Real's fourth generation of "ancestor doubling" collectively dubbed Sedenions, seemed to become fatally overrun with "monstrosities" – was paralleled, in work contemporary with Hurwitz's (and the ultimate source of Lévi-Strauss's), by the beginnings of de Saussure's structural analysis of linguistic signs. This quickly led to Louis Hjelmslev's "Net" of indefinitely many strata, each housing a bifurcation into "form" and "substance" dichotomies (etic *vs* emic, as in "phonetic" *vs* "phonemic"), later given a name in its own right by André Martinet: the indefinitely extensible process of "double articulation."[6]

## 1. Hjelmslev's Net and Box-Kite Basics

The semantic sieving by binaries called "Hjelmslev's Net," "whose minimal number" of successive distinctions "guarantees their maximal extension," starts with a split into identically defined, arbitrarily named categories of *expression line* and *content line* which "have mutual solidarity through the sign function."[7] Deleuze and Guattari depict this cascade of "reciprocal presuppositions" (Roman Jakobson's term for it) this way:

> Now this net had the advantage of breaking with the form-content duality, since there was a form of content no less than a form of expression. Hjelmslev's enemies saw this merely as a way of rebaptizing the discredited notions of the signified and signifier, but something quite different was actually going on. Despite what Hjelmslev himself may have said, the net is not linguistic in scope or origin (the same must be said of double articulation: if language has a specificity of its own, as it most certainly does, that specificity consists neither in double articulation nor in Hjelmslev's net, which are general characteristics of the strata).[8]

From here, the "structural semantics" of Algirdas Greimas' "Semiotic Square" will show us how to link Catastrophe Theory to the foundations of semantics via ZD "representations," with an assist from the researches of Jean Petitot, disciple of Greimas and René Thom both. But for those familiar with Deleuze's other writings, something of the sort is already implicit in the quote above: for the term "strata" is in fact inspired by Catastrophists' special use of it, and as such was the focus of an earlier work on Leibniz, named after the simplest Catastrophe, The Fold. (And cited at length in numerous places germane to our interests, by Petitot.[9] Recall, too, it was Leibniz who adapted the ancient Chinese scheme of "solid" and "broken" lines, in the $2^6$ hexagrams of the I Ching, seeding all later ideas of binary bit-strings!) This is made explicit in the authors' introduction to "double articulation":

> God is a Lobster, or a double pincer, a double bind. Not only do strata come in at least pairs, but in a different way each stratum is double (it itself has several layers). Each stratum exhibits phenomena constitutive of *double articulation*. Articulate twice, B-A, BA. This is not at all to say that the strata speak or are language based. Double articulation is so extremely variable that we cannot begin with a general model, only a relatively simple case. The first articulation chooses or deducts, from unstable particle-flows, metastable molecular or quasi-molecular units (*substances*) upon which it imposes a statistical order of connections and successions (*forms*). The second articulation establishes functional, compact, stable structures (*forms*), and constructs the molar compounds in which these structures are simultaneously actualized (*substances*). In a geological stratum, for example, the first articulation is the process of "sedimentation," which deposits units of cyclic sediment according to a statistical order: flysch, with its succession of sandstone and schist. The second articulation is the "folding" that sets up a stable functional structure and effects the passage from sediment to sedimentary rock.[10]

But how might ZD's relate to Catastrophes? Consider, first, in Catastrophe Theory, the simplest "compact" unfolding of a Singularity (meaning, no "tune in next time" sendings to or from "infinity"– i.e., beyond the model's "local scoping") entails the splitting of a stable regime, or "attractor," into two. These may compete until capture of one by the other, and exhibit paradoxical interplay *en route* (e.g., the "girl playing hard to get" surrendering when her suitor gives up, or the "rat backed into a corner" jumping the cat).[11] But simply reverse the sign on the potential function which defines the "behavior" associated with twiddling this Cusp's two "controls," and we get a "Dual Cusp," with two unstable *maxima* ("mountain peaks" or "fountains," not "basins of attraction"), mediated by a solitary zone of stability (instead of a compromise-inhibiting "zone of inaccessibility," necessitating the standard Cusp's trademark "jump"). This is the dynamics of a *pan-balance*, which seesaws until the weights on its two arms are in equilibrium (a "measurement" is obtained). Similarly, when two non-zero quantities have a Zero product, this result can be interpreted in the manner of "pair production" from the far-from-empty QM Void. In ZD's generated by CDP, a zero product of 2 of them isn't just "nothing," but a pair of oppositely signed copies of a 3rd ZD, which always zero-divides the first 2's "diagonal twins" (see next paragraph). If we use "~" to write "mutually zero-divides," we can define an *equivalence relation*, peculiar to ZD's, upon which we'll construct an NKS substitute for "Pythagorean triplets" (actually, "six-packs") using letters A – F.

To preface a longer story with some shorthand, the letters stand for vertices of special kinds of triangles, called Sails. Each letter further denotes a plane ("Assessor") spanned by two imaginary axes, whose diagonals are saturated with (and in fact contain *all* the "primitive") ZD's. 4 such Sails, joined only at 6 shared vertices, comprise half the faces of an Octahedral vertex figure, or "Box-Kite" – a 3-D iconic shorthand for the 12-D space ZD's *really* span! Its 4 "empty" triangles are Vents where "$2^N$-ion winds" blow. If we arbitrarily label one Sail's vertices A, B, C, and those of its opposite Vent F, E, D respectively, then the 3 orthogonal lines AF, BE and CD demarcate "Struts": stabilizing dowels made of wood or plastic in real-world Box-Kites.[12] For all and only the vertex pairs at opposite ends of a Strut, mutual ZD's don't exist. Along all 12 edges of the Box-Kite, however, the "~" relation obtains between vertices, and indicates a binary selection rule: if the edge is marked with a plus sign ("+") on the Box-Kite, then similarly oriented diagonals at each vertex – both passing through the same 2 quadrants of their respective planes ("/" for the {– –, ++} or "synchronized"[13] line; "\" for the "anti-synchronized" or {– +, + –} line) – are mutual ZD's . If marked with a minus sign ("–"), however, the reverse obtains: a slash ("/") at either vertex mutually ZD's the other's backslash ("\").

It is an invariant feature of Box-Kites that one Sail will have all 3 edges marked "–", with its opposite Vent likewise having 3 "–" edge-signs, while all 6 remaining edges are marked "+". By convention, the all-minus Sail, dubbed the "Zigzag," is labeled A, B, C, making the all-minus Vent F, E, D. The Zigzag gets its name from contemplating the sequence of diagonals joined in ZD pairings as one circuits its vertices: if we start with a "/," then we get a 6-cyclic "double covering" of this sort: "/ \ / \ / \". The other 3 Sails, with 2 "+" edges and 1 "–", do a double-circuit that suggests the under-over tracings of the knots they're named for: Trefoils ("/ / / \ \ \"). The "3 + 1" logic of Trefoils *vs* Zigzag will be seen to play a fundamental role in semantic models. To understand it, though, we must now show where the XOR-ing of integer indices comes into play between edge- and strut- joined vertices, and within the pairs of axial units forming planes at each.

From an infinity of regular forms in low dimension (e.g., "polygons"), only triangles and squares generalize to N-space (ergo, simplexes and "Cartesian coordinates"). Likewise, while there are many ways to represent "multiplication tables" of hypercomplex numbers (2 for Quaternions: $ij = k$ or $ik = j$; 480 for Octonions, labeled via 7 nodes and 7 arrows on the PSL(2,7) triangle; 100's of billions for Sedenions), only *one* scheme works consistently, without anomalies, for *all* CDP $2^N$-ions. Index all $2^N$ units by integers, from 0 (for "reals") to $2^N - 1$ (hence, up to 15 for Sedenions). The *product* of 2 such will have as *its* index the XOR of the indices of the units *producing* it. A Sedenion "for instance" is found in the triplet (7, 12, 11): that is, $i_7 * i_{12} = i_{[7 \text{ XOR } 12]} = i_{11}$.

*Signing* of products is more involved; in fact, the anomalies it inevitably engenders are what make ZD's possible! Signing is simple in one crucial instance: "adding a bit to the left" to $2^N$-ions' indices effectively multiplies them "from the right," by the unit indexed $2^N$ – thereby generating $2^{N+1}$-ions. $2^N$ hence is the "generator" (henceforth, G) for the $2^{N+1}$-ions; within the latter, index XORing tells us that for $k < G$, $i_k * i_G = + i_{(k + G)}$, and oppositely for $k > G$. For all other indices, associative triplets (i.e., clones of Quaternions, or "Q-copies," with 7 in the Octonions, 35 in the Sedenions, and $(2^N - 1)(2^N - 2)/6$ for general $2^N$-ions) display conflicts between "signing" and "counting" order. This is due to the "NKS effect" of iteratively adding bits via successive "G actions," N growing. In an XOR-compliant writing of Octonion triplets, if the first 2 indices are in the order that generates a "+" sign for the 3rd, then (1, 2, 3) – "*i, j, k*" indices for standard Quaternions – have "sign ordering" identical to "counting order" as shown; ditto, for (1, 4, 5); (2, 4, 6); (2, 5, 7); and (3, 4, 7). But the 2 remaining triplets, if written to show "sign ordering," must necessarily violate "counting order": (1, 7, 6) and (3, 6, 5). And while it is possible to fiddle with labeling schemes to reduce the number of such "out-of-sync" triplets to 1, one can never expunge them entirely: there is always a bump in the rug. And, starting with the Sedenions (viz., the (7, 12, 11) triplet cited last paragraph), such bumps will clash in such a manner as to mandate the existence of Zero-Divisors.

Corollarily, such ZD's will *never* involve products including the perfectly behaved G. This leads to a surprisingly simple rule for finding all the "primitive" ZD's in the Sedenions: take one unit whose index, L, is less than G, then take another whose index U exceeds it, *and does not equal the XOR of G with L* (equivalent to saying, in this special case, it doesn't equal G + L). For any of the 7 choices of L, there are then 6 choices left for U, such that the plane they span (one of the 7 * 6 = 42 "Assessors") will have diagonals of the form $k * (i_L \pm i_U)$, $k$ real, all of whose points are ZD's. Put another way, point coordinates on these lines (hence, $k$ values) are irrelevant, ergo arbitrary: ZD's are therefore best thought of as *line elements* – meaning they, like NKS entities generally, as Wolfram has argued the case, have an inherent randomness *built into their workings*.[14]

This rule (with restrictions due to "carry-bit overflow") generalizes for higher $2^N$-ions; the point of moment, for present purposes, is that Sedenions partition into 7 Box-Kites, each housing 6 Assessors, each of which forms ZD's only within the 2 Sails sharing its vertex in its Box-Kite. This is the concrete manifestation of Moreno's abstract "$G_2$ homomorphism": since $G_2$ symmetry is either 14-fold (when the subscripted "2" – for the count of unit vectors spanning the space containing its dozen "roots" – is included, as when particle physicists calculate symmetries) or 12-fold (as in coding theory or kaleidoscope building), the simplest homomorphism to contain *both* possibilities has an 84 count – the number, precisely, of ZD-saturated lines, when tallying across all 7 Box-Kites.[15]

## 2. "Greimas's Glitch" and the ZD-Net Hypothesis

Although this may seem like reaching on a first reading, the last few paragraphs set up the second resonance between ZD's and the fundaments of semantics. The first, of course, was the suggestive analogy between the nonstop redoubling of Hjelmslev's Net and the Cayley-Dickson Process. The second, though, is readily made by those who know what to look for: for the dyadic nature of ZD's – each comprised of two indecomposable "units," one with index less than the generator G, one with index greater, which is a formal way of saying they inhere in different analytic "levels"– is just what Algirdas Greimas tells us is called for to resolve the most deep-set contradiction in the classical structuralist agenda. "Greimas's Glitch," as I like to think of it, can be depicted thus.

Once upon a time, there was a fundamental epistemological problem in Saussurean signification theory. It was celebrated in all its naïveté on the cover of the Winter '72 issue of the postmodernist review Diacritics, which displayed a picture of the wrapper from a then-new candy from Peter Paul, a peanut-butter bar called "No Jelly." This was meant to be juxtaposed with a lift from a text of near-scriptural status, Ferdinand de Saussure's Course in General Linguistics, printed on the inner back cover: "In language there are only differences. Even more important: a difference generally implies positive terms between which the difference is set up; but in language there are only differences *without positive terms.*" (This number became famous for other reasons, too: it contained the first major English interview with deconstructionist Jacques Derrida – one of whose principal targets, of course, was the ensemble of de Saussure's unexamined assumptions!)

The credo of pure arbitrariness in signs (any random signifier *could*, in theory, be linked to any random signified) meant their patterns of mutual difference *alone* made signifying thinkable – which made the existence of synonyms and their generalizations a *problem*. Greimas used the "double articulation" agenda initiated by Hjelmslev to posit a way out of this embarrassment: assume two "pieces" of a signifying complex reside at different strata, so that they essentially could work together when combined with another such complex to yield a "zero" (hence, synonym effect) in semantic difference. (Think non-message bits in a self-correcting error code.) This point is too far-reaching *not* to instantiate it: in Structural Semantics, Greimas gives as example the "synonymous" morphological marks "-s" and "-en," for designating plurals ("cows" *vs* "oxen") in English. These different marks can cover identical elements of content, "on the condition, however, of having different contextual distributions."

> It is consequently sufficient that a mark be realized twice, at two different structural levels – the first time under the form of the opposition of phonemes, the second time, under the form of the opposition of phonematic segments – for the differentiating effect of the first mark to be canceled by the advent of the second differential gap. Thus in some conditions two marks can be neutralized by combining with each other:
>
> $$X + (-X) = 0$$
>
> and redundant variation on the plane of expression does not provoke any difference of signification. Thus synonymy is possible.[16]

Thus, whatever other significances we may attribute to it, Zero here indicates a state of "no difference" between semantic contents, and hence a redundancy which context may frame as lack or ambiguity of meaning, but which *change* of context may convert into a *bifurcation* of possibly intended senses, as with puns, jokes, and ironic statements – with the affect responses these can induce being minimally modeled by standard Cusps (whose nonlinear "jumps" in such cases, between the attractors competing to be "taken seriously" by us, are frequently flagged by sudden bursts of laughter).[17]

Primitively speaking, potential meaning is linked with potential for affect release – the "Aha!" experience, that is, contains the gist of "making sense." But we also know two other things. First, "level-jumping" is indigenous not just to specific "meaning experiences," but has a volatility independent of any such instances, which allows us to navigate between "candidates" for meaningfulness with a freedom verging on infinite: in a slogan, *anything can be quoted or pointed at ("**That!**")*. Such pronominal "pointing" serves precisely to fulfill its own name: any phenomenon, however complex and hence needful of enormous numbers of dimensions to be modeled, can be magically shrunk to a point – with the reverse mountain-from-molehill process, of course, obtaining (*"**point well taken!**"*) just as naturally. Derrida calls this "iterability," noting that "Every sign, linguistic or nonlinguistic" is amenable to citing/sighting/siting/inciting us to insight: "This citationality, duplication, or duplicity, this iterability of the mark," he tells us, "is not an accident or anomaly, but is that (normal/abnormal) without which a mark could no longer even have a so-called 'normal' functioning."[18]

What this means for our ambitions here is that our modeling must have built-in "open-ended-ness," with the indefinite extensibility of $2^N$-ions as its vehicle. We will ultimately need *means* for collapsing high-order ensembles into "points" – and not just means, but *motives*: there must, that is, be mechanisms indigenous to our models which enable "explosions" of insight-flashes into epics. (Poincaré puts his foot on the bus, and sees in a flash how his Fuchsian groups' transformations can "provide transportation" for hyperbolic geometry as well!)[19] Before we can even consider what that might require, however, let's note that its fulfillment must contradict the basic assumptions which buttress up our choice of where to start from – which leads me to the second "other thing" referenced in the paragraph before the citation. $2^N$-ions offer us structural alternatives extending far beyond the limits of "Bott periodicity"; but the Catastrophist's toolkit is "local" in the extreme, with the "elementary" theory being heavily constrained by the 8-D "Bott limit": "new tricks" stop with the 8-D end of the 3-term E series of singularities.

A famous conundrum, in fact, with implications for all finite mathematical classifications, had its origins in this very domain. I refer to Vladimir Arnol'd's 1972 "A, D, E Problem,"[20] which showed the taxonomy of local singularities Catastrophes provide in fact complies, at a meta-level, with countless other schemes, shorthanded by the "ball-and-stick" drawings physicists call "Dynkin Diagrams," and algebraists interested in finite reflection groups (a.k.a. *kaleidoscopes*) know by Coxeter's name. In this setting, the trivial accelerated object in a gravitational potential, the simple Fold of sudden starts and ends, the doubly complex (but still quite simple) Cusp, are just the three lowest forms on the infinite "totem pole" of "$A_n$" type, with n = 1, 2, 3 respectively. But as the guy painting the street meridian put it in "Chafed Elbows," "You have to draw the line someplace" – and this will be where *we* draw *ours* to start with. (Merging this "starter kit" into NT's fractal "Big Picture" will concern us later – and indeed, is a key target.)

Let's start, then, by noting that if Zero-Divisors can indeed serve as markers to indicate a *lack of difference* (or, dually, an *ambiguous "balancing act"*) between two meanings, in what manner can they facilitate a *basis* for locking a meaning in place? Corollarily, if simple Catastrophe models can underwrite simple "meaning effects," in what sense do ZD's and these models intersect?  As a first pass at an answer, let's posit the following five-step "ZD-Net Hypothesis":

1. Since ZD's seem, in some primary way, to facilitate synonymy (and hence *inhibit* individuation of significance), the most fundamental "meaning effects" should devolve upon just those substructures of ZD ensembles – "Sails" and "Box-Kites" and (starting with the 32-D Pathions) multiple-Box-Kite synergies – wherein ZD's are *not* acting "as such."

2. The most basic ZD objects which *don't* act as such are "Strut-opposites":  only *oppositely* placed Assessors on a Box-Kite do *not* mutually zero-divide, and their 4 units obey an "XOR arithmetic" equivalent to the logic of Greimas' "Semiotic Square" (SS).  As Petitot makes the case that the SS, a fundamental "crossing" of two semantic polarities, can be approached as the simplest "stacking" interaction of *two* simple 2-Control Cusps (which is to say, in terms of the fifth term and its dual in the A series, *one* simple 4-Control *Butterfly*), this implies the standard and dual $A_5$ can be "represented" by "Strut logic," with the horizontal, diagonal and vertical paths on the SS corresponding to its invariant features (Strut Constant "S," Generator "G", and their XOR "X", not necessarily in that order).

3. Sails contain substructures which are *not* zero-divisors:  the "L" units (indices < G) of each Assessor in a Sail form a Q-copy; but each also forms a Q-copy with the "U" units (indices > G) of its two "Sailing partners."  For the Zigzag ABC, we can write these this way:  {a, b, c}; {a, B, C}; {A, b, C}; {A, B, c}.  The same obtains for Trefoils, but the subtle "Trip Sync" property says that "bow-tie flows" (akin to Greimas' "figure-8's") through their shared Assessor, from one Zigzag Q-copy to another in a Trefoil, allow for a special kind of D-series "message passing."  This process also permits an interpretation of a second SS reading, favored by American postmodernists, and (due to his claiming an isomorphism between his Square and the Klein Group) sanctioned by Greimas himself.

4. The eliciting of Tests and Helpers and their interrelationships are the hallmark of Greimas' analysis of Vladimir Propp's study of Russian wondertales (the study from which the SS construct emerged).  As Greimas points out, there are 3 of each in any such tales' architecture, and this can be related to the 3 Struts (and 3 "Trip Sync" message-pass setups, one per Trefoil) in a Box-Kite.  The "contract," meanwhile, which enables the tale as a whole, can be associated with the ZD-free Q-copy with units indexed by S, G, X and 0 (the real unit) – what a C++ or Java programmer would call an "abstract class," which the tale's unfolding instantiates.

5. The simplest CT object which incorporates all the above is also the simplest "non-elementary" Catastrophic form – one rigorously related by Arnol'd's school to the Equilateral among the 3 "primordial triangles" in Plato's <u>Timaeus</u>.[21]  (The other 2 underwrite "higher harmonics" of Box-Kite representation, thanks to which the pathway to NT embeddings via higher $2^N$-ion cascading is opened to our inquiry.)

### 3. Semiotic Squares and "Eco Echo"

The Italian writer Umberto Eco is best known for his novels – one of which, <u>The Name of the Rose</u>, was made into a highly successful Sean Connery vehicle. But he is also a major European academic, and one of the world's pre-eminent semioticians. In his best-known scholarly text, he's noted a fundamental difficulty in the enterprise of nailing down the "basics" of what meaning is, which will be referenced herein as "Eco echo":

> as soon as the semantic universals are reduced in number and made more comprehensive … they become unable to mark the difference between different sememes. And as soon as their number is augmented and their capacity for individuation grows, they become *ad hoc* definitions…. The real problem is that *every semantic unit used in order to analyze a sememe is in its turn a sememe to be analyzed.*[22]

The problem is rendered far less paradoxical by seeing semantics as side-effect of a mathematical substrate (one of the primary motivations for representing fields of knowledge as mathematical "sciences" in the first place!). Yet "Eco echo" still has a residuum of "Cheshire cat's smile" adhering to it, which must be somehow addressed in the mathematical apparatus we put in place. Complementary to Derrida's "iterability," "Eco echo-ing" requires that we at least treat the 3 distinct binary categories of relationship which constitute the Semiotic Square as having a protean nature, which double articulation can underwrite (and find, within this process, its own impetus). Relations represented by SS verticals (perhaps strut-constant "S" values in Box-Kite talk) can, in other SS's serving other purposes, appear as horizontals (hence, G xor S = G + S = "X" values), or even as diagonals (generators themselves, once "relabeling schemes" are allowed) in others.

Allied with such conversions are fundamental patterns of "explosion," the simplest of which show how G, S, and X patterns among Sedenion Box-Kites can be mapped to ensemble effects in the 32-D Pathions' various septets (7 "Pléiades" plus the "Atlas" or "Muster Master") and triads (the 7 "sand mandalas") of Box-Kites – which we'll express as special "blow-ups" along S, G and X pathways, with each having its own somewhat surprising properties. As the aim here is to be illustrative rather than exhaustive, we'll be providing exemplifications of patterns amenable to greater generalization, rather than anything like an "axiomatics." (We will later, when we speak of higher-order "semical harmonics" of Box-Kite representational analysis, consider ways in which SS's themselves can be seen as "lines" within "meta-SS" diagrams – but again, our aim will be to introduce, not exhaust, such topics, which are after all still in their infancy as research themes, yet should be taken as candidates for indefinite induction anyway.)

This much said, some "buffer refresh," via a new Box-Kite diagram and ancillary tables of index patterns, is next, with attention directed initially to peculiarities obtaining between Strut-Opposite dyads (as well as between the diagonals within a single Assessor at either end of such a dyad). In the diagram, uppercase demarcates indices > G; lowercase, those less. A, B, C mark the Zigzag (on top), and D, E, F the Vent opposite it, with the negative edge-signs framing them drawn in red ink to make a Star of David (with the outer hexagon of positive edge-signs in blue). Zigzag ABC, and Trefoils ADE, FDB and FCE, are painted moon-silver, blood-red, sun-yellow and sky-blue, respectively.

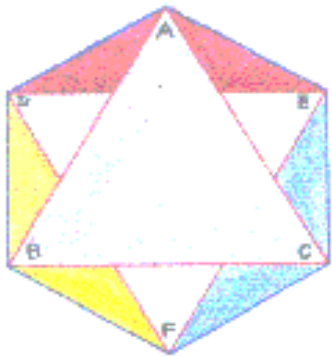

| Strut | Assessors at Sedenion Box-Kite Vertices | | | | | |
|-------|------|------|------|------|------|------|
| **Const** | **A** | **B** | **C** | **D** | **E** | **F** |
| **1** | 3,10 | 6,15 | 5,12 | 4,13 | 7,14 | 2,11 |
| **2** | 1,11 | 7,13 | 6,12 | 4,14 | 5,15 | 3, 9 |
| **3** | 2, 9 | 5,14 | 7,12 | 4,15 | 6,13 | 1,10 |
| **4** | 1,13 | 2,14 | 3,15 | 7,11 | 6,10 | 5, 9 |
| **5** | 2,15 | 4, 9 | 6,11 | 3,14 | 1,12 | 7,10 |
| **6** | 3,13 | 4,10 | 7, 9 | 1,15 | 2,12 | 5,11 |
| **7** | 1,14 | 4,11 | 5,10 | 2,13 | 3,12 | 6, 9 |

Note that, in this way of rendering the Box-Kite, the products of the lowercase letters joining vertices along the outer blue hexagon have sign-orders that are always positive when tracing a *counter*clockwise path: for strut constant 1, for instance, a*d = 3*4 = + 7, the e at the ADE Trefoil's 3rd vertex. But now observe that the reverse is the case for products of UPPERCASE letters along the same paths: A*D = 10*13 = –7, the e with *reversed* sign this time (e and not E, since the high bits are XOR'd away). For the edges in red along the Star of David, abc *and* ABC products are positive along counterclockwise pathways, while def and DEF are positive along *clockwise*. (See BK4 for a brief explanation and proof.) This is the basis of the "Trip Sync" property alluded to above: "slippage" is allowed from one Q-copy to the other sharing a lowercase index *if it belongs to the Zigzag* (equivalent to saying, as we'll be seeing, that exchanges along an X path in a Semiotic Square are "free"). Otherwise, a "quantum of spin" (or orientation reversal) will result. In the uniquely "all red-edged" Zigzag Sail, such slippage is *always unobservable* in this context. (The importance of this will become clear once we translate Strut Opposites' multiplicative patterns into SS language. "Slip spin" reverses relations 2 Squares share, in a way involving D-series Catastrophes – with "meaning effects" involving Greimas' trio of "Helpers," enabling one SS "hand" to "wash the other.")

A "/" or "\" times itself, at any vertex, is some scalar times the negative real unit. The product of "/" and "\," in the same Assessor plane, is ± a scalar times X ( = G + S). Next in complexity are Strut-Opposite products: writing (V, v) and (Z, z) for "def" Vent and Zigzag Assessors which straddle a Strut, the rules are just V * Z = v * z = S; Z * v = V * z = G; Z * z = v * V = X. But to which of the Semiotic Square's 4 corners (yielding 4! = 24 choices) should each go? And which arrow-types should be labeled G, S and X? (6 more choices!) Might the answers to these questions depend, to some extent, upon the specifics you wish to elicit? It would do to describe the 3 path types in classic Jakobson-Greimas language, then try out a couple of mappings while watching what happens . . . but first, let's come at things from the other side of our overriding comparative urge, and derive the Square's own context by an iterative process – CDP-driven "double articulation" of the Reals – then (per Wolfram's mantra) "watch what happens."

To inaugurate that last exercise, an initial setup should suggest itself: the indefinitely extensible morass of ZD's, caught up in Flows beyond the confines of metrics, can be thought of as the "sediment" of the first Deleuzian articulation; the Lie-algebraic "elementary taxonomies" conforming to "A-D-E" diagrams provide the Forms, or "sedimentary rock" of the second. We could call the latter "belonging structures," and the former "meaningful flows" – carriers of the shapes and colors Hjelmslev's Net goes fishing for.

If you've been reading the discursive footnotes (especially 15, and 19-21), or are familiar with the literature on "A-D-E" already, then you know that this divide between Forms and Flows is hardly arbitrary. All Forms are positive-curvature bound, which is how the McKay Correspondence is literally able, by projection from SU(2)'s 4-space to the 3-D sphere Plato himself knew, to map polygons to the A-series, dihedra (polygons colored on both sides) to the D, and the E-series to the 3 symmetries underlying the 5 regular solids (with cubes and octahedra, and icosa- and dodeca- hedra, being "duals" of each other, such that the vertex figure of either maps to the faces of the other). The "primordial triangles" of interest to us, per the fifth part of the ZD-Net Hypothesis, form a thin "boundary layer" of "flat" or "Euclidean" curvature, containing – in a dimension-reversing sense topologists should feel at home with – the 3 Platonic Solids as maximal strata: the equilateral triangle, the tetrahedron ($E_6$); the right isosceles, the Box-Kite's octahedron ($E_7$); the 30-60-90, the icosahedron (with "closest packing" of 8D spheres mapped by Octonionic loops of 240 "integral Cayley numbers," the $E_8$ "root count").

These are also the cross-sections of the 3 prismatic kaleidoscopes one can build with 3 mirrors (one for each of the "Timaeus" singularities' behavior variables – and, in particular, one for each Strut in our Box-Kite). These singularities are all "non-elementary," which is to say they have "parameters" which measure the breakdown of the link between topology and differentiability – which mysterious variables, based upon projective invariants (e.g., perspective's cross-ratio for the non-Equilaterals), are rendered palpable and even prosaic by prismatic scoping: instead of a single center of symmetry, as in the standard 6-fold kid's symmetry-maker, you get honeycombs tessellated by the triangle in question, with as many centers as your viewing aperture will let you take a peek at. This "modularity" (related, but not identical, to the notion of that name we talked about at the outset) is linked to certain critical values which, when hit upon, cause "explosions" into infinite-dimensional "control spaces" – just the sort of impetus our discussion of "iterability," at the start of last section, said we'd be needing.

But needing for what, exactly? For throwing our Net into the negative-curvature spaces of higher-dimension reaches, where Chaotic "period-doubling" and other modes of turbulence hold forth. The non-elementary figures of what Form-like patterns exist there are all related, as Dolgachev showed, to M. C. Escher's world of hyperbolic tilings – and, as Jacques Hadamard first pointed out in the year of Hurwitz's Proof, such spaces of negative curvature are generally *divergent* (in senses that Poincaré's then-ongoing efforts were just making discussable).[23] Add to this, the seemingly gratuitous fact that James Callahan, in his now-classic multi-dimensional navigations of $E_6$ and Double Cusp Catastrophes using stacked 2- and 3-D "tableaus," has found marvelous defining features he calls "trance tunnels" and "ship's prows," through and with which one "steers" in "control space," by use of "keels" and other sailing instruments.[24] Now imagine our flat boundary layer of "Timaeus triangle" unfoldings to be akin to the surface of deep waters plied by fishing trawlers – with the fishermen throwing out seines we've come to call by Hjelmslev's monicker. Then imagine the sea across which we're sailing to be a "sentient ocean," like that covering the surface of Stanislav Lem's <u>Solaris</u>. In the sci-fi classic, this realm of Meaningful Flows can respond to our most secret wishes and thoughts with instantly realizable architectural (and even sensual) manifestations. Something of the sort is what we hope will be in store for us, as we begin to cast our net for Semiotic Squares and other "deep-sea creatures," in higher-dimensional Zero-Divisor spaces.

## 4. Struts, Tray-Racks, and Twisted Sisters: Bounding the Model from Above

Let's start with the usual "assume true" of so many mathematical proofs. Assume it's the case, that is, that a Strut is a showcase for Semiotic Squares, as claimed. Let's look simple-mindedly at what that means geometrically. On the vertex net we've called a Box-Kite, the 3 mutually orthogonal Struts are themselves diagonal lines – 2 at a time in each of the 3 mutually orthogonal Squares whose perimeters contain all 12 Box-Kite edges. So we immediately see that the context in which we claim the Semiotic Square should be housed has another level of "double articulation" already enabled. These doubled-up SS's, or "Tray-Racks" – named for those portable, snap-together-fast make-shift table-tops, so perfect for holding TV dinners while your attention's focused on something else – can have their perimeters traced in the manner of Sails. But aside from the obvious 4- vs. 3- fold symmetry, there's at least one other very big difference: since the edge-signs regularly alternate from plus to minus as you circuit, there is no double covering; instead, one has *two* distinct circuits which never intersect, each "making zero" with each of the 4 Assessors through which it threads. But where one Tray-Rack does so with the slash, the other does its drive-by at the same site with the backslash.

There is also a double articulation of product outcomes: in Sails, the $3^{rd}$ vertex has "pair production" of its contents result from multiplication of the other two. With Tray-Racks, each edge belongs to a different Sail, but a sort of "pair production" still results – with a "twist." In fact, in BK1, I called this second kind of multiplication just that: a "twist product." The ways in which it "double articulates" are best explained by defining it: if $(i_U + i_L)(i_{U*} + i_{L*}) = 0$, then $(i_U + i_{L*})(i_{U*} - i_L)$ does, too. (I leave off the tedious variations based on internal signing.) What we have, then, is a swapping rule: a same-case-indexed unit of one partner in a mutual ZD pair is swapped for the corresponding unit in the other, with the resultant ZD pairings having edge-signs toggled. (More, just as in SS's, horizontal and vertical edges have different behaviors, with parallels of one type "twisting" in tandem to opposite-type parallels on different Tray-Racks.) The edge in question, meanwhile, is in a *different Box-Kite* – although, as we'll clarify next, the Sail containing this alien edge is in the same *Loop* as the Sail we started with.

The divide between Loop and Box-Kite, in fact, is the basis of this "double articulation" effect. Loops are just the non-associative algebraist's next best thing to Groups, and Raoul Cawagas, after but independently of me, also eked out the patterns of ZD's in the Sedenions, via the Loop-studying software he'd designed for just such purposes.[25] Except his 7 hitherto undiscovered 16-element Loops (ZD-harboring clones of Octonions, which – being uninterested in Loop theory at the time of BK1 – I'd simply given the silly name of "automorphemes," when I wasn't just referring to them as "Moreno counterfeits") partition into 4 subsystems, corresponding to Sails. These Loops, then, and my Box-Kites, are duals: each Box-Kite's 4 Sails inhabit different Cawagas 16-Loops, and vice-versa. Moreover, the patterns governing this duality are exceedingly orderly, and can be catalogued with PSL(2,7) triangles I call "Twisted Sisters" diagrams, whose 7 nodes house the Strut Constants of the Box-Kites engaged in "swappings," while edges have color-coded segments (one per each Tray-Rack) to flag in- and out- going matchups.

As the phenomena just limned serve to bound our SS modeling "from above," let's stop to take a look at what's just been sketched out so hurriedly. (Later arguments about higher "semical harmonics" will require, after all, some familiarity with such matters.) First, let's look at twist products, from one Tray-Rack to each of two others.

To see the orderliness of the Twist, a standard presentation of the Tray-Rack is needed. Regardless of Box-Kite, all 3 Tray-Racks share a structural invariance: unlike the 4 Sails, which display a "3+1" symmetry between Trefoils and Zigzag in their cycling patterns, each Tray-Rack displays a "3+1" in the orientation of ZD flows. Adopting an "overtone series" convention that low bits are less volatile, hence less subject to movement, we use them as points of reference.[26] In each Tray-Rack, the orientation of L-bit products is counterclockwise along 3 successive sides, with the 4th (with negative edge-sign) showing a clockwise reversal. If we label them with the letters of the struts perpendicular to them in the octahedron, Tray-Racks AF, BE, CD have reversed edges DE, FD, EF respectively, each of which belongs to a Trefoil completed by the perpendicular's Zigzag Assessor. (The "3+1" of Sails and Tray-Racks hence share the same source.) Rotate their frames so that the top edge and the projection of the perpendicular on the square's center define this Trefoil, and paint it in traffic-light colors: red, yellow, or green for Tray-Racks labeled AF, BE, or CD. Call such presentations "Royal Hunt" diagrams, after the 5th line text of I Ching Hexagram 8, "Holding Together": "In the hunt, the king uses beaters on three sides only and foregoes game that runs off in front."

The vertical lines, if twisted together (a "V*" operation), produce a Tray-Rack in another Box-Kite, and likewise for the horizontals ("H*"). Further, performing H* or V* twice is isomorphic to the workings of an associative triplet, indexed by the 3 Strut Constants of the Tray-Racks' Box-Kites. Since twist results can be written in 2 different orders, depending on whether the L or U units are swapped left to right, a convention to minimize rewriting and maximize consistency is adopted: for top and left edges in "Royal Hunt" setup, L's are swapped and U's held in place; for bottom and right edges, U's are swapped and L's fixed. When this procedure is followed, H* produces a Tray-Rack rotated 90° clockwise of "Royal Hunt," while V* produces one rotated an equal but opposite amount. Depicted below are 2 of the 3 twist product sets for the S=3 Box-Kite.

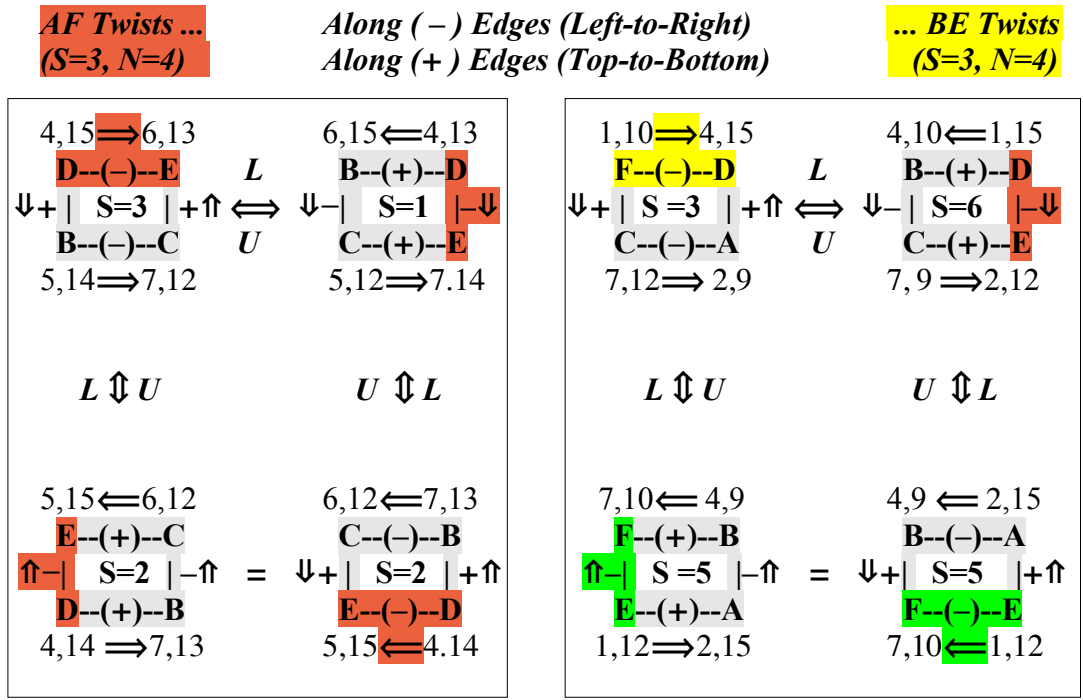

In general, the S of the Box-Kite we twist *to*, is easily shown to always equal the index of the perpendicular's Vent Assessor, in the Tray-Rack twisted *from*. In the particular cases just shown, the S-triple (1,2,3) for Box-Kite 3 indicates all Tray-Racks are of the same AF variety: this is unique to the "true" Quaternion index set, and explains the central red circle in the "Twisted Sisters" diagram. The "Royal Hunt" sequence for the same Box-Kite's BE Tray-Rack, by contrast, requires all 3 colors for its (3, 6, 5), as do most S-triples. Since Strut Constants for $2^N$-ions are always less than $2^{N-1}$, the Sedenions' "Twisted Sisters" diagram has '4' (i.e., the Octonions' 'G') in the center, where it plays a role opposite the one it adopts with index multiplication: as *symmetry breaker* (viz., 347 and 246's 2-coloring, 145's color ordering). All other paths are quite regular. We've seen the "inner circle" of Quaternion indices forms an all-red circle of clockwise flows. The nodes at medians point inwards with green; yellows point into the vertices. The clockwise flowing on the perimeter always leaves a node with green, passes through yellow, then comes to a stop at its terminus with red (the "traffic light" sequence).

Such PSL(2,7)-based diagrams exist in the Pathions as well, shown in BK4 by the "Stereo Fano" synergies of paired PSL(2,7)'s shorthanding the layout of the 14-Assessor, 7-Box-Kite ensembles (7 Pléiades plus the "Atlas"), and can serve to refine the classifying of Box-Kites, now that the Strut Constant proves no longer unique, hence too coarse.

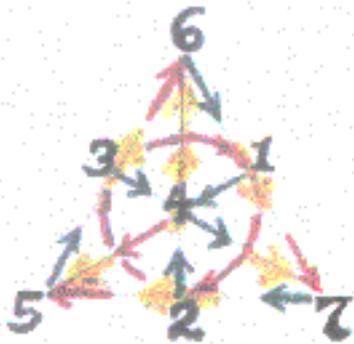

We will not be concerned with them here, but we *will* care about the general phenomenon of studying structures like Sails and Tray-Racks which thread ZD's, and which we'll refer to in the abstract as "lanyards" – those strands of leather, string or plastic upon which charms or beads are queued, to make closed loops of jewelry. They roughly do service in ZD work comparable to what Groups, built on '1' and notions of Identity, do in standard algebra and number theory. Informally, lanyards are linked lists such that each diagonal of each Assessor threaded is left-multiplied to make zero once only, with closure.

The Sail is a doubly linked list, with all diagonals in the 3 Assessors visited being threaded in one 6-cycle; Tray-Racks, also doubly linked, only provide single covering, hence come in pairs, each cover threading half the diagonals without intersecting the other.

Two more lanyards are highly significant: the Quincunx threads 5 Assessors, and can be interpreted as a strategy to go around the obstacle to traffic flow the top edge of a "Royal Hunt" diagram represents: by detouring up the Strut perpendicular to the Tray-Rack, then coming back to the endpoint of the bypassed segment, a 10-cycle double covering can be achieved without "going against the flow" of the natural "L-bit" orientation. A doubly astonishing connection exists between the symmetry of choices available and the icosahedral reflection group: *doubly* so, as the connection of the "obstacle bypass" problem in Singularity Theory to this same group caused a big to-do when members of Arnol'd's circle discovered it (then heard, as well, that diagrams of the wavefront indicating this symmetry had appeared in the first-ever analysis text by l'Hôpital, who likely derived them ultimately from Huyghens himself!).[27] Finally, all 12 diagonals are threaded by the polyrhythmic "Bicycle Chain" of ¾-Tray-Rack scans, linked by "–"-edge jumps: an inviting candidate for NKS-style simulation work – about which, more toward the end.

## 5. First Pass at the Model, and the Route to Fractals

The correspondence of Semiotic Squares as modeled by ZD subsystems with the CT treatment of Petitot is further reinforced by the fact that the "ceiling" of our approach, as indicated by the just-discussed Tray-Rack-based "swapping" problem twist products handle, is *also* at a level above *his* SS model. As an index of the limits of CT modeling itself, one must employ the highest one-behavior stratum of the Double Cusp, the $A_7$ "Star," to effect such bartering – and even then, the result is not quite satisfying![28]

This shared "bounding from above" of both ZD-based and Petitot models is underscored by the Tray-Rack's relationship to the second (right isosceles) "Timaeus Triangle," and hence the Double Cusp, whose "extended Dynkin diagram" is a box whose 4 corners house "basins of attraction" separated by saddles on their edges, and joined by diagonals which cross in a mountainous "maximum" (and whose so-called "generic section" in Control Space is the 4-cusped curvilinear square known as an "Astroid").[29] But do we truly need 16-D to contextualize both? We'll want to consider how the models are jointly "bounded from below." But to do so, we need an initial presentation of our own, in Greimas's terms.

To begin at the beginning: in Greimas's <u>Structural Semantics</u>, whose analysis of Propp's wondertale morphology led directly to the SS diagram, Greimas asks us what, on the linguistic plane, we mean when we say we "perceive differences"? He answers his own question in a doubly articulated fashion:

> 1. To perceive differences means to grasp at least two object-terms as simultaneously present.
> 2. To perceive differences means to grasp the relationship between the terms, to link them together somehow.[30]

We are next told that "simultaneous grasping" implies its own dyad: first, one object alone doesn't carry signification; second, this perceiving together implies both something shared between them (the "problem of resemblance") and (in order to discern there are, indeed, two or more of them) something somehow different, which serves to distinguish them. Resemblance and difference, then, are part of what is simultaneously present – within perception *as such* – *before* we can speak at all of *what* it is we've perceived. (Not to mention, the *how* and *why* of it.) We then must make our starting point the notion of *reciprocal presupposition* (Hjelmslev's *solidarity*), as within the percept as such, the one object-term can carry no signification without the other. This is not only the most important assumption Greimas requires of us; it "is the only undefined logical concept which allows us to define reciprocally, following Hjelmslev, signifier and signified."[31]

But Greimas and Hjelmslev both want to transplant this notion to virtually every stratum of description: the theory of signification, then, is but a "surface effect"; what they're after would use this sole "undefined logical concept" in its arsenal as an all-purpose tool. To reach the Semiotic Square from here, one "double articulates" what we've just been given, and assumes that *meaning* is the result of the interlinking of two such dyads, giving us the Square's four corners. But the verticals, which effect the soldering between the initially unrelated dyadic "semes," do *not,* by dint of their opportunistic nature, entail *reciprocal* presupposing, but only the one-way notion of "implication," with the lower term somehow "containing" or being "generatively linked to" the one above it.

Greimas also tells us that the diagonals, which indicate "contradiction" between a seme and the absence of same, is linked to another of Jakobson's 3 primary categories of binary opposition, called "privation": 2 phonemes, otherwise identical, differ in the presence of a mark: 'g' is voiced; 'k' is not. Is not the simplest way to indicate this by a simple on/off bit – what "G toggling" gives? Giving a commonsensical yes to this question quickly tells us our S should mark the vertical "implication" pathway. Here's why.

First, mapping G to the 'privative' pathway means, once we've made our way through the blizzard of names Greimas and his interpreters have used to highlight the various attributes each pathway is assumed to display – as a relation between "paradigmatic" place-holders or "states" on the semantic plane; as "syntagmatic" segments in a performance of state-to-state syntax, itself dependent upon an "actantial" analysis which further discriminates "actants" from "roles"; etc. – that we see it simply as the pure form of the foundational Box-Kite binary, *mutual exclusivity*, implicit in saying 2 terms contradict each other: if Castor is present, Pollux is not. But G alone, for reasons already covered, offers the possibility of linking this XOR notion to a no-exceptions, "+" or "–".[32]

Next, *reciprocal presupposition* is the basis of linking the two semes into one semantic bundling, and this operation is inherently (and explicitly described this way by Greimas as) *hierarchical*: conjunction on one level, disjunction on the other. Horizontal pathways are thus the *only* places that can say "X" marks the spot: for we then have the ZD *diagonals themselves* on one level, and the *units spanning their plane* (and defining their equations as just $k(i_U \pm i_L)$) on the other. By elimination, "S" must be the vertical.

But one can argue positively, as well: given a seme, the first thing we frame is a notion of its absence (hence, the "contradiction" pathway: mice only play in the lower right slot if the seme "cat" is absent in the upper left; if a "light" seme is switched off, barracks-bound recruits start sleeping – or maybe vampires leave their coffins). So, before we can consider the possibility of fusing two semes into a Square (which "mergers and acquisitions" work is what the "implication" process underwrites), we already have "G" in hand – and, for any given CDP level (which is what "G," as a power of 2, suffices to indicate) the necessary precondition for such union is "S." For Sedenions, all "S" values generate successful candidates for ZD status when our simple production rule is used; for Pathions and higher $2^N$-ions, however, carry-bit overflow puts this status in question, if S exceeds 8 and is not equal to some higher power of 2. For Pathions, this means all S > 8 yield ensembles of 3 Box-Kites instead of 7, with the 7 distinct patterns having a clear "flip-book" sequencing determined strictly by simple unit in- or de- crements to S. For the 64-D Chingons, 3 distinct ranges display different-sized ensembles. As N increases further, we essentially are faced with a "spectroscopy" problem.

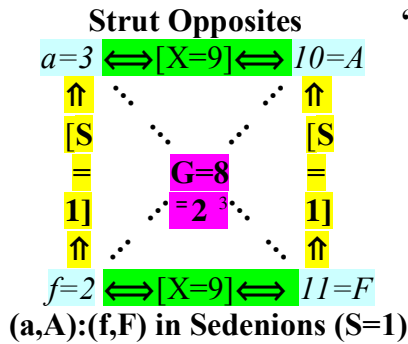

**Strut Opposites**

$a=3 \Longleftrightarrow [X=9] \Longleftrightarrow 10=A$

[S = 1]  [S = 1]

G=8 =2

$f=2 \Longleftrightarrow [X=9] \Longleftrightarrow 11=F$

**(a,A):(f,F) in Sedenions (S=1)**

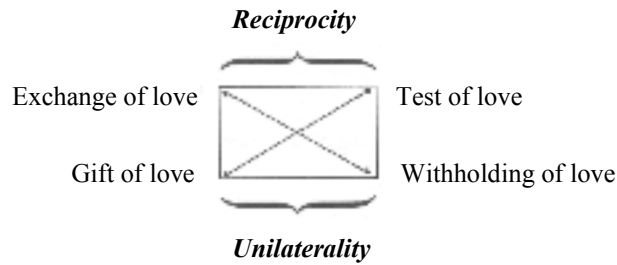

**"FiduciaryLove"** *à la* **Racine's Roxane in** <u>Bajazet</u>
(A Semiotic Square from Greimas & Fontanille[33])

*Reciprocity*

Exchange of love — Test of love

Gift of love — Withholding of love

*Unilaterality*

To make this seem concrete, and see how it relates to attaining the NT "fractal limit," a review of "emanation tables" (on which, see BK2) is in order. These do for ZD's what Cayley Tables do for group operations. For a particular stratum of $2^N$-ions, K = $(2^{N-1} - 2)$ Assessors interact in Box-Kite ensembles for a given Strut Constant value: 6, 14, 30, 62 for N = 4, 5, 6, 7, and so on. For any S value, regardless of carry-bit overflow, we can make a K x K product table, place Assessor labels on each row and column head, and enter the results in the appropriate cells. Define headings in left-right, top-down order, entering lowest L-index first, and put that of its strut-opposite in the mirror-opposite slot; put the next-lowest not yet entered to the right of the first (and its strut-opposite just left of the first such); continue in this manner until all L-indices are accounted for.

If Assessors indicated by headings for a cell mutually zero-divide, put their *emanation* (instead of "0") in the cell; otherwise, leave it empty. For all such tables, both long diagonals will form an empty "X", since strut-opposites and an Assessor's own diagonals never mutually zero-divide. For S values not subject to overflow, all other $K^2 - 2K$ cells ( = 24, 168, 840 for N = 4, 5, 6) will be filled. As this number is always a multiple of 24 (the number of cell entries defining a Box-Kite), this means Strut Constants, beyond the Sedenions, underwrite *ensembles* of interconnected Box-Kites: 7, 35, 155, and 651 for N = 5, 6, 7, 8 (also the number of associative triplets for N = 3, 4, 5, 6).

For S values where carry-bits matter, things get more complicated. As explained in BK2, the 7 Pathion emanation tables for 8 < S < 16 are more sparse than expected: only 72 cells (= 3 Box-Kites) are filled. With Chingons, 3 broad classes emerge in 64-D, the first of which – due to its obvious "fractality" – we'll examine in detail later on.

8 < S < 16 gives 19 instead of the "full set" of 35. Within this range, take the Sedenion "sand mandala" with the same S, chop it into quadrants, and attach copies of each to the starter image's corners. For S=15, fill the 8 same-sized regions on the perimeter with mirrorings of the starter-image quadrants, but with the empty diagonal stretches replaced by 16's (or, if the diagonal below contains 8's not 7's, by "strut-opposite" 31's), and all other indices *augmented* by 16. Analogous "period-doubling" strategies obtain in higher reaches, indicating an emanation-table-based route to the fractal limit as N increases. (See drawings on left, below, for S=15 in Pathions, then Chingons.[34]) For 16 < S ≤ 24, we've 7 Box-Kites, whose 14 x 14 central square mimes the Pathion emanation table for S – 8, *sans* the crossing pairs of near-solid parallels – making the central square for S=24 *empty*. 24 < S < 32 gives 23 Box-Kites (for S=25, see below right – and note how "white space," unlike for S=15, is *filled with*, not *empty of*, entries!).

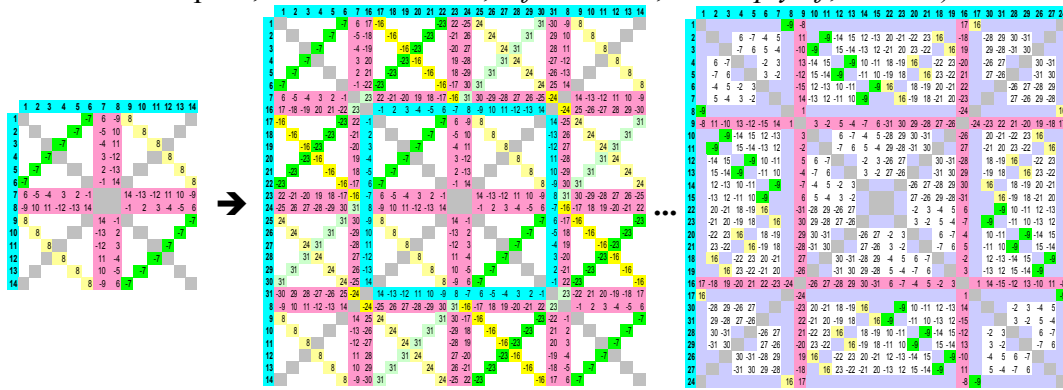

**S=15: 32-D Pathions ➔ 64-D Chingons ( ... and, S=25, among the Chingons only)**

**Appendix: Box-Kite "Do-It-Yourself Kit"**

The basic "Box-Kite Machine" is an interactive form with 2 slots and 2 buttons: enter the power of 2 which generates the $2^N$-ions of interest, then enter the Strut Constant (any value from 1 to $2^N - 1$) of the Box-Kite ensemble you wish to explore. Click the "Produce Tone Row" button, and the lower parts of the form fill in in a split second: low and high indices of all the Assessor pairs, in emanation-table row and column, with the full spelling out of these dyads in the "Tone Row Pairs" listing at the end. After you've spent about two seconds getting that much produced (with maybe 200 lines of code doing the work "under the hood"), you click one more button to "Create Emanation Table" according to the specifications. Except for N pushing two digits or worse, the CSV text-file dumped out by this operation takes a split second, and then you just roll over into an Excel spreadsheet, import said textfile into a worksheet, and then play with colorizing.

Ignoring user-interface overhead ("exiting event" code in field slots and other stuff that's part of the "implementation cost" but otherwise irrelevant), the whole machine runs on only 3 little programs (ca. 100 lines each, written in the Visual Basic customized for the Notes object-hierarchy called "LotusScript" that the business analyst's workstation I was building for a major investment house required I immerse myself in). I've attached the highly inelegant code for these below. (Inelegant, because written superfast, under duress of having things I was supposed to be doing that paid the bills; also, because the library calls were subject to some parameter-passing bugginess that I had to do grotesque workarounds to bypass, leading to very stupid-looking code. The point, though, was to get this machine built in a hurry, and not worry about "elegance"!)

This could easily be ported to Java or a Mathematica notebook, where much more sophisticated things become possible (e.g., New Kinds of Databases and simulation work based on them). The 3 modules do the following:

1) Crunch Cayley-Dickson Process into a stupidly self-evident method for generating signed products of arbitrary imaginaries of any subscript (see the discussion of the 3 rules that run everything in my prior "BK4" paper on the Wolfram Science site: "Triple Negatives Rule" and "Trident Product" and "Recursive Bottoming" and that's it).

2) Given N and S, create the appropriate lower-index terms to serve as row and column labels in the "emanation table," and then put them in proper order (being sure, of course, to eliminate the "strut-constant" S itself from said listing!).

3) A double-for-loop which runs through the row-and-column listing of indices twice, calculating the cell values in the table. There's 2 arrays for low- and high- index units in the Assessors each such label signifies, and then the full Assessor of each row label is multiplied times each column's full Assessor, and if they DO zero-divide, the third Assessor in their Sail has its lower-index term entered in the appropriate cell. (If the two Assessors do NOT mutually zero-divide, the cell is left blank).

That's it! And if Part II's proof of the all-important Sky meta-fractal's existence and nature seems too skimpy, just look at the code-produced results I'll showcase there. Pascal's notion of "mathematical induction" will make "proof" a given.

The first function performs signed multiplications of two hypercomplex units, given their subscripts (i.e., the "XOR of indices = index of product" rule is assumed). Here's the code, which can be thought of as "Cayley-Dickson for Dummies":

**Function M( LI%, RI% ) As Integer**

```
Dim QSigns( 0 To 3, 0 To 3) As Integer
Dim NegTally%, LBits%, RBits%, G%
QSigns( 0, 0 ) = +1
QSigns( 3, 1 ) = +1
QSigns( 1, 3 ) = - 1
For i% = 1 To 3  ' 0 => index of real unit; save for 0 * 0 = ( + R), all index-0
        ' products of +-signed units = ( - R)
        QSigns( 0, i% ) =  +1
        QSigns( i%, 0 ) =  +1
        QSigns( i%, i% ) = -1
        If i% < 3 Then
                QSigns( i%, i% + 1 ) = +1
                QSigns( i% + 1, i% ) = - 1
        End If
Next
NegTally = +1
If LI% < 0 Then ' Keep it simple!  After these if-tests, can assume we're dealing
        ' with absolute values of indices.
        NegTally% = ( -1 ) * NegTally%
        LI% = - LI%
End If
If RI% < 0 Then
        NegTally% = ( -1 ) * NegTally%
        RI% = - RI%
End If
XorRoot% = LI% Xor RI% ' Except for all-important signing, this is the index of
                ' the product of the units whose indices are LI% and RI%
Do While True  ' Time for recursion
        LBits% = Len( Bin$( LI% ) )
        RBits% = Len( Bin$( RI% ) ) ' Number of bits - and hence, smallest power
                ' of 2 which exceeds - the integer in question
        If LI% = 0 Or RI% = 0 Then
                Exit Do
        End If
        If LI% = RI% Then
                NegTally% = ( -1 ) * NegTally%
                Exit Do
        End If
        If LBits% < 3 And RBits% < 3 Then ' We've got a product of Quaternion
                ' indices:  extract sign from QSigns array and we're done!
                NegTally% = QSigns( LI%, RI% ) * NegTally%
```

```
                    Exit Do
            End If
            If LBits% = RBits% Then  ' Two indices arise from same
                    ' generator = 2^(LBits% - 1) = G
                    G% = 2^(LBits% - 1)
                    If LI% = G% Then
                            Exit Do ' LI% < RI%, and triplet =
                                    ' ( LI% Xor RI%, LI%[=G%], RI%)
                    End If
                    If RI% = G% Then
                            NegTally% = ( -1 ) * NegTally%
                            Exit Do ' Reverse of prior if's situation
                    End If
                    If ( LI% Xor RI% ) = G% Then
                            If RI% > LI% Then
                                    NegTally% = ( -1 ) * NegTally%
                            End If
                            Exit Do ' Triplet always has form (lo, G, hi)
                    End If
                    NegTally% = ( -1 ) * NegTally%
                    LI% = LI%  - G%
                    RI% = RI% - G% ' For generator G, row-index = LI% = G + r
                    ' and col-index = RI% = G+c, r & c < G, LI%*RI% = (-1)*r*c
                    Goto RECURSIVE
            End If
            If LBits% < RBits% Then
                    G% = 2^(RBits% - 1)
                    If RI% = G% Then
                            Exit Do
                    End If
                    If ( LI% Xor RI% ) = G% Then
                            NegTally% = ( -1 ) * NegTally%
                            Exit Do
                    End If
                    NegTally% = ( -1 ) * NegTally%
                    RI% = RI% - G% ' LI% < G% and so is left unchanged
                    Goto RECURSIVE
            End If
            If RBits% < LBits% Then
                    G% = 2^(LBits% - 1)
                    If ( LI%  Xor  RI% ) = G% Then  ' e.g., 15 * 7 = 8 => (7, 8, 15)
                            Exit Do
                    End If
                    NegTally% = (-1) * NegTally%
                    If LI% = G% Then ' e.g., 16 * 3 = ( -19 ) => (3, 16, 19)
                            Exit Do
```

```
                        End If
                        LI% = LI% - G% ' reverse logic of last if-test
                        Goto RECURSIVE
                End If
RECURSIVE:
        Loop
        M = NegTally%  *  XorRoot%
End Function
*****************************
```

Emanation Table creation is handled by the next piece of code:  it lives in the "click event" of the button you can see in the screenshot, "Generate Tone Row."  The name is inspired by Schonberg's tone-row, and the choice of strut-constant-based row of Assessor labels does service as the "musical scale" from which all else is "composed."

**Sub Click(Source As Button)**

```
        Dim s As New NotesSession
        Dim w As New NotesUIWorkspace
        Dim ui As NotesUIDocument
        Dim doc As NotesDocument
        Set ui = w.CurrentDocument
        Set doc = ui.Document
        Dim  N_test$, S_test$, powmax%, generator%, rowmax%,
                sconst%, xconst%, i%, count%, opp%
        N_test$ = Trim$( ui.FieldGetText( "Nion" ) )
        If N_test$ = "" Then
                powmax% = 5 ' Default = work with Pathions
        Else
                powmax% = Cint( N_test$ )
        End If
        generator% = Cint( 2 ^ (powmax% - 1) ) ' For 32-D Pathions, = 16
        S_test$ = Trim$( ui.FieldGetText( "Strut" ) )
        If S_test$ = "" Then
                sconst% = 1
        Else
                sconst% = Cint( S_test$ )
        End If
        xconst% = generator% + sconst%

        rowmax% = Cint( generator%  - 1 ) ' For 32-D Pathions, = 2^4 - 1 = 15
                           ' = highest "LO" index
        Redim lstRaw%( 1 To ( rowmax% - 1 ) ) ' Low-index tone row elements, in
                           ' counting, not "strut-opps at opposite poles," order
        Redim lstLoTone%( 1 To ( rowmax% - 1 ) ) ' Eliminate strut constant ('S')
                           ' from the count
        Redim lstHiTone%( 1 To ( rowmax% - 1 ) ) ' Eliminate XOR of strut constant
                           ' with generator ('X') from count
```

```
                count% = 0
                For i% = 1 To rowmax%
                        If i% = sconst% Then
                                Goto JUMP
                        End If
                        count% = count% + 1
                        lstRaw%( count% ) = i%
JUMP:
        Next
        Dim lo_count%, hi_count%, try%
        lo_count% = 1
        hi_count% = rowmax% - 1
        For i% = 1 To ( rowmax% - 1 )
                try% = lstRaw%( i% )
                x% = ( try%  Xor  sconst%   )
                If try% < x% Then
                        lstLoTone%( lo_count% ) = try%
                        lstHiTone%( lo_count% ) = ( try%  Xor xconst% )
                        lstLoTone%( hi_count% ) = x%
                        lstHiTone%( hi_count% ) = ( x%  Xor xconst% )
                        If 2 * lo_count% = ( rowmax% - 1 ) Then
                                Exit For
                        End If
                        lo_count% = lo_count% + 1
                        hi_count% = hi_count% - 1
                Else
                        Goto SKIP
                End If
SKIP:
        Next
        doc.ToneRowLo = lstLoTone%
        doc.ToneRowHi  = lstHiTone%
End Sub
```

******************************

Finally, once you've picked G and S, and generated the Tone Row, you want to generate a CSV-file of raw ASCII text, from which importing to Excel will dump out a table of even gargantuan proportions without much trouble – and then you can "make it pretty" with Excel's easy cell-format commands.

**Sub Click(Source As Button)**
```
        Dim s As New NotesSession
        Dim w As New NotesUIWorkspace
        Dim ui As NotesUIDocument
        Dim doc As NotesDocument
        Dim recBuf As String
```

```
Dim prod, Nion, HI_NDX As Integer
Dim H%, K%, L%, Q% , fnum%, idx%
Dim LRow%, HRow%, LCol%, HCol%
Dim varLo( ) As Integer
Dim varHi( ) As Integer
Set ui = w.CurrentDocument
Set doc = ui.Document
Nion = doc.Nion( 0 )
If Nion < 10 Then
        strN$ = "00" & Cstr(Nion)
Elseif Nion < 100 Then
        strN$ = "0" & Cstr(Nion)
Else
        strN$ = Cstr(Nion)
End If
Strut = doc.STRUT( 0 )
If Strut < 10 Then
        strS$ = "00" & Cstr(Strut)
Elseif Strut < 100 Then
        strS$ = "0" & Cstr(Strut)
Else
        strS$ = Cstr(Strut)
End If
FileNm$ = "C:\Documents and Settings\Bob\My Documents\Box-Kites\N" &
        strN$ & "S" & strS$ & ".txt"
fnum% = Freefile( )
Open FileNm$ For Output As fnum%
HI_NDX = ( 2^(Nion - 1) - 2 ) ' For Sedenions, have 8 - 2 = 6 cells per
                ' emanation-table row or col; for 2^5-ions, 16 - 2 = 14;
                ' for 2^6-ions, 2^5 - 2 = 30; ...
Redim Preserve varLo( 1 To HI_NDX ) As Integer
Redim Preserve varHi( 1 To HI_NDX ) As Integer
idx% = 0
Forall lo In doc.ToneRowLo
        idx% = idx% + 1
        varLo( idx% ) = lo
End Forall
idx% = 0
Forall hi In doc.ToneRowHi
        idx% = idx% + 1
        VarHi( idx% ) = hi
End Forall
recBuf = ""
For L% = 1 To HI_NDX
        recBuf = recBuf & "," & Cstr( varLo( L% ) )     ' Top row of col-labels
Next
```

```
For k% = 1 To HI_NDX
        LRow% = varLo( k% )
        HRow% = varHi( k% )
        recBuf = Cstr( LRow% )
        Dim Edge%, Edge2%  ' if +, then edge-sign => / * / OR \ * \ for the ZD
                        ' current; if -, then / * \ OR * / \ /
        For q% = 1 To  HI_NDX
                If ( q% = k% ) Or ( q% = HI_NDX + 1 - k% ) Then
                        recBuf = recBuf & ","
                        Goto ZDNEXT
                End If
                LCol% = varLo( q% )
                HCol% = varHi( q% )
                ' ZD_Test -- modularize into a library function later!!
                firstnum% = HRow%  ' will this workaround do the trick?
                (Internal LotusScript parameter-pass library-call bugginess!!)
                nextnum% = LCol%
                UpperLeft% =   M( firstnum%, nextnum%)
                firstnum% = HRow%
                nextnum% = HCol%
                UpperRight% = M( firstnum%, nextnum% )
                firstnum% = LRow%
                nextnum% = LCol%
                LowerLeft% =   M( firstnum%, nextnum% )
                firstnum% = LRow%
                nextnum% = HCol%
                LowerRight% = M( firstnum%, nextnum% )
                If Abs( UpperLeft% ) = Abs( LowerRight% )  Then
                        If Sgn( UpperLeft% ) = Sgn( LowerRight% ) Then
                                Edge% = ( +1 )
                        Else
                                Edge% = ( -1 )
                        End If
                Else
        ' Something's screwy!  These should ALWAYS have same abs. values?!
                        Print "Product of " & Cstr( UpperLeft% ) & " and "
                        & Cstr( LowerRight% ) & "doesn't work right?!"
                        recBuf = recBuf & ", ??"
                        Goto ZDNEXT
                End If
                If Abs( UpperRight% ) = Abs( LowerLeft% )  Then
                        If Sgn( UpperRight% ) = Sgn( LowerLeft% ) Then
                                Edge2% = ( +1 )
                        Else
                                Edge2% = ( -1 )
                        End If
```

```
                    Else
            ' Something's screwy!  These should ALWAYS have same abs. values?!
                            Print "Product of " & Cstr( UpperRight% ) & " and "
                            & Cstr( LowerLeft% ) & "doesn't work right?!"
                            recBuf = recBuf & ", ??"
                            Goto ZDNEXT
                    End If
                    If Edge% = Edge2% Then
                            ' We have a pair of mutual zero-divisors!
                            recBuf = recBuf & "," & Cstr( Edge% * Ab( LowerLeft%))
                            ' Print low-index, with edge-sign, of emanation
                    Else
                            recBuf = recBuf & ","
                    End If
ZDNEXT:
            Next
            Print #fnum%, recBuf
        Next
        Close #fnum%
        Msgbox "Text file " & FileNm$ & " successfully created.", 64, "EMANATION-
TABLE TEXT-FILE"
        End Sub
```

1    Lifted from http://www-groups.dcs.st-and.ac.uk/~history/Quotations/Laplace.html

2    The title of a highly influential precursor of "small-world" theory, Mark S. Granovetter's paper in The American Journal of Sociology, 78 (1973), 1360-1380, whose original context was, in fact, the way people find jobs:  not, as one might at first assume, through their closest friends (since the overlap between their and one's own connections will be too great to provide too many "leads"); but rather, through remote acquaintances (who have their own small groups of close friends, with which one is likely to have almost no intersection *except* through such "weak ties").  "Small-world" theory as such, meanwhile, derives from a paper whose influence has proven inversely proportional to its length:  Duncan J. Watts and Steven H. Strogatz, "Collective dynamics of 'small-world' networks," Nature, 393 (1998), 440-442.  For a fuller presentation, see Duncan J. Watts, Small Worlds:  The Dynamics of Networks Between Order and Randomness (Princeton NJ:  Princeton University Press, 1999); for the intimate relationship between small-world and scale-free network theory, see Albert-László Barabási, Linked:  How Everything Is Connected to Everything Else And What It Means for Business, Science, and Everyday Life (New York:  Penguin Books, 2003).  "Six degrees of Kevin Bacon"– a take on the title of John Guare's 1990 similarly-themed play, Six Degrees of Separation – derives from three inebriated fraternity brothers who managed to ride their goofy idea into TV and Internet exposure:  see Ann Oldenburg, "A thousand links to Kevin Bacon:  Game calculates actor's connection," USA Today (October 18, 1996), p. 51.  One of the clearest instances of a "small-world network," ironically, is provided by the Kevin-Bacon-like "Erdös number" of mathematical collaboration:  those who've written papers with the prolific Hungarian have Erdös number = 1; those who haven't collaborated with Erdös himself, but have co-authored papers with people who have, have  Erdös number = 2; no mathematician who's ever collaborated on a paper has been found with  Erdös number higher than 17, and most – over 100,000 – have  Erdös number 5 or 6:  see Jerry Grossman's " Erdös Number Project Home Page," www.oakland.edu/~grossman/erdoshp.html.  This is ironic, because Erdös collaborated on a famous series of papers that created the theory of random networks – a beautiful theory which "small-world" thinking had to chuck out the window before it could find itself eliciting laws for *real* networks.  Those with an obsession should also visit Brett Tjaden's "Oracle of Bacon" web page, which contains a list of the 1000 best-connected actors and their distribution sequences (many of the most connected of whom, surprise, are porn stars):  www.cs.virginia.edu/bct7m/bacon.html.

3    A. R. Rajwade, Squares (Cambridge UK:  Cambridge University Press, 1993), especially Chapter 10:  "The $(r,s,n)$-identities and the theorem of Hurwitz-Radon (1922-3)," 125-40, and the background provided in Chapter 1:  "The theorem of Hurwitz (1898) on the 2,4,8-identities," 1-11.  For further background, see I. L. Kantor and A. S. Solodovnikov, Hypercomplex Numbers:  An Elementary Introduction to Algebras (New York:  Springer-Verlag, 1989), an elegant introduction to the whole theme of hypercomplex numbers, all building up to the formulation, in Chapter 17, of Hurwitz's 1898 proof that "Every normed algebra with an identity is isomorphic to one of the following four algebras:  the real numbers, the complex numbers, the quaternions, and the Cayley numbers."

4    R. Guillermo Moreno, "The zero divisors of the Cayley-Dickson algebras over the real numbers," Bol. Soc. Mat. Mex. (3) 4:1 (1998), 13-28; http://arXiv.org.abs/q-alg/9710013.

5    Claude Lévi-Strauss, The Naked Man:  Introduction to a Science of Mythology/Volume 4, John and Doreen Weightman, transl. (New York:  Harper & Row, 1981; French original, 1971), pp. 644-5:  "It would seem that the point at which music and mythology began to appear as reversed images of each other coincided with the invention of the fugue, that is, a form of composition which, as I have shown on several occasions, exists in a fully developed form in the myths, from which music might at any time have borrowed it.  If we ask what was peculiar about the period when music discovered the fugue, the answer is that it corresponded to the beginning of the modern age, when the forms of mythic thought were losing ground in the face of the new scientific knowledge, and were giving way to fresh modes of literary expression.  With the invention of the fugue and other subsequent forms of composition, music took over the structures of mythic thought at a time when the literary narrative, in changing from myth to the novel, was ridding itself of these structures.  ***It was necessary, then, for myth as such to die for its form to escape from it, like the soul leaving the body, and to seek a means of reincarnation in music.***"  [Bold italics provided by this author, not the anthropologist.]

6    André Martinet, Elements of a General Theory of Linguistics, L. R. Palmer, transl. (Chicago:  University of Chicago Press, 1964; French original, 1960), p. 24 et. seq.  A useful précis can be found in an accessible postmodernist text that can provide a good introduction to the general issues and vocabulary we'll be working with, Ronald Schleifer's A. J. Greimas and the Nature of Meaning:  Linguistics, Semiotics and Discourse Theory (Lincoln NB:  University of Nebraska Press, 1987).  On p. 17, he explains that "'Double articulation,' a term developed by André Martinet, is the presence – the reciprocal presup­position – of both combination and contrast. In this conception, language presents two distinct planes of analysis. Each of the 'units of the first articulation,' Martinet argues, 'presents a meaning and a vocal (or phonic) form. It cannot be analyzed into smaller successive units endowed with meaning . . . But the vocal form itself is analysable into a series of units each of which makes its contribution to distin­guishing *tête* from other units such as *bête*, *tante*, or *terre*.' In traditional linguistics, the first articulation is the combinatory of grammar (or morphology) while the second is that of phonology." As with CDP, redoubling can be extended indefinitely – and, in fact ("Hjelmslev's Net"), *must* be.

Louis Hjelmslev, Prolegomena to a Theory of Language, Francis J. Whitfield, transl. (Madison WI: University of Wisconsin Press, 1963; Danish original, 1943), pp. 59-60.

[8] Gilles Deleuze and Félix Guattari, A Thousand Plateaus: Capitalism and Schizophrenia, Brian Massumi, transl. (Minneapolis and London: University of Minnesota Press, 1987; French original, 1980), p. 43. The authors further "topologize" Martinet and Hjelmslev's basic insights: "One is never without the other: a double deterritorialization..." (P. 392) And, in percipient conformity with the current theme, they develop a philosophy of Number they call "nomadic" (cutely opposed to Leibniz's "monadic") which is primordially topological: "It is the number itself that moves through smooth space.... The more independent space is from a metrics, the more independent the number is from space." (P. 389). Or again, "Smooth space is a field without conduits or channels. A field, a heterogeneous smooth space, is wedded to a very particular type of multiplicity: nonmetric, acentered [i.e., not coordinated by Identity operators], rhizomatic multiplicities [tracing connections between inherent and palpable Singularities, rather than branchings from a unified Tree] that occupy space without 'counting' it and can 'be explored only by legwork.' [e.g., as embedded in a *network*]." (P. 371)

[9] Of most moment, for our interests, is Deleuze's short essay, "A quoi reconnait-on le structuralisme?," in F. Chatelet, ed., Histoire de la philosophie (Paris: Hachette, 1973), which Petitot discusses at length in Section 3.3, "Deleuze's proposal for a schematism of structure," of his Morphogenesis of Meaning, transl. Franson Manjali (New York: Peter Lang, 2004; French original, 1985), henceforth MM. The notion of "place-holder" is argued to be fundamental to all structuralism, and the "empty place" the basis of any notion of "symbolic." As Petitot notes, referring to this piece, "The identity of a symbolic place is not what ensures its stability, but what ensures the possibility of its displacement.... If the rela­tive displacement (metonymy) can be an intrinsic part of the identity of position, it is because every structure 'contains an object or an element which is quite paradoxical.' This paradoxical element is of a kind very different from the symbolic elements, the differential relations, and the singularities. It circu­lates within the series as if it was 'its own metaphor and its own metonymy' (p. 322). It lacks any onto­logical function (it is not an object), any self-likeness (it is not an image), and logical identity (it is not a concept) (p. 323). And if the relative displacements is [sic] an intrinsic part of positional identities, it is because the relative place of the terms in the structure depends on their absolute place in relation to this element." (Petitot's p. 61). But what can this 'quite paradoxical ' element be but Zero in a place-based scheme of digitization, with structuralism thereby being an exercise, ultimately, in ZD representings?

[10] Ibid., pp. 40-41. Geological imagery is also one of the core motifs of René Thom's CT magnum opus, Structural Stability and Morphogenesis (Reading MA: W. A. Benjamin, 1975; French original, 1972), henceforth SSM: of the 27 illustrations in the 1975 translation's "centerfold," 6 are of such geological "catastrophes" as fault-lines, beach domes, river beds, and superposed laminar strata in rock; many more such images populate the main body of the text (viz., the drawings accompanying the 3[rd] endnote to Chapter 7, showing "threshold stabilization" at work in the erosion of a jurassic anticline).

[11] A wealth of homey "toy models" of this variety can be found in E. C. Zeeman's now-classic Scientific American article of April, 1978 – especially the 60-page draft version included in his Catastrophe Theory: Selected Papers, 1972-1977 (Reading MA: Addison-Wesley, 1977), pp. 1-64, as well as the "Social Sciences" and "Discussion" chapters (302-407, and 604-658, respectively).

[12] For what it's worth, this is a literal fact: octahedral lattices are known to be rigid; an octahedral box-kite with pairs of dowels joining along each strut axis in the center, effectively produces an 18-edged figure with 7 vertices. Per A. L. Loeb's degrees-of-freedom formula, rigidity requires $E \geq 3V - 6$, which means we get just the right count even if we somehow assume 3 dowels can magically (i.e., in purely "mathematical space") pass through the central 7[th] vertex without requiring their splitting into pairs of real-world sections. See the discussion on rigidity in Jay Kappraff, Connections: The Geometric Bridge Between Art and Science (New York: McGraw-Hill, 1991), pp. 270-3.

[13] The term is not chosen arbitrarily: anyone who's used an oscilloscope knows that the appearance of diagonals indicates the two orthogonal variables represented onscreen are "in sync"; it is also known

that synchronization plays a highly organizing role in Chaotic systems. On all this, see Steven Strogatz, <u>Sync: The Emerging Science of Spontaneous Order</u> (Hyperion Books, 2003), especially Chapter 7: "Synchronized Chaos."

[14]  This point is developed further in my third "Box-Kite" monograph (henceforth, "BK3"), "Quizzical Quaternions, Mock Octonions, and Other Zero-Divisor-Suppressing 'Sleeper Cell' Structures in the Sedenions and $2^N$-ions," available at http://arXiv.org/abs/math.RA/0403113. For background on the basics of Box-Kites, index labeling schemes, and CDP, see http://arXiv.org/abs/math.GM/0011260 ("BK1"), "The 42 Assessors and the Box-Kites they fly: Diagonal Axis-Pair Systems of Zero-Divisors in the Sedenions' 16 Dimensions." The logic and graphics of "carry-bit overflow" in higher $2^N$-ions' "emanation tables" – specifically, the "sand mandalas" in the 32-D Pathions, which provide the basis for extrapolations to scale-free networks – are treated in "BK2": "Flying Higher Than a Box-Kite: Kite-Chain Middens, Sand Mandalas, and Zero-Divisor Patterns in $2^N$-ions Beyond the Sedenions," http://arXiv.org/abs/math.RA/0207003. Finally, a general introduction, and harnessing of ZD thinking to NKS concerns, can be found in "BK4": my NKS 2004 presentation, "The 'Something From Nothing' Insertion Point: Where NKS Research into Physics Foundations Can Expect the 'Most Bang for the Buck'," http://www.wolframscience.com/conference/2004/presentations/materials/rdemarrais.pdf.

[15]  Gell-Mann's "Eightfold Way," based on SU(3) symmetry – equivalently, $A_2$ in "Dynkinese," which yields up the simplest Fold Catastrophe with germ $X^3$ and one control working on the X monomial – gets its count by taking the hexagon of roots (hence, the correspondence to the most primitive kid's kaleidoscope, made by 2 mirrors placed at a $60^o$ angle, in "Coxeter" dialect) and adding 2. If the 2-D backdrop is assumed, closest packing of spheres within it (the link to coding theory) is just the root-count, 6: for $A_3$, the 3-D closest-packing number of 12 is the root-count as well (or the image-count in a 3-mirror kaleidoscope, mirrors 1 and 2, 2 and 3 meeting in dihedral angles of $60^o$ to make a wedge). For 4 and 5 dimensions, though, the closest packings correspond to the D-series root-counts associated with so-called "Feynman checkerboard" patterns, of 24 and 40 respectively (28 and 45, if "Cartan subalgebras" are included). For 6-, 7-, and 8- D, the E-series serves similarly, with respective root-counts of 72, 128, and 240 (78, 133, and 248 with subscripts ladled in). Superstring theory's first great "aha!" occurred when the "anomaly count" for 16-D heterotic strings was shown to be 496: just the tally generated by two orthogonal spaces, each associated with the Octonions' (non-associative) algebra – or "$E_8$ x $E_8$" for short! Meanwhile, on the theme of "84" which triggered this note: it is perhaps ironic that this number is associated with the study of "Hurwitz Groups," which play a crucial role in the study of compact Riemann surfaces of genus at least 2 (whose automorphism groups have maximum order of 84 (g–1)). It's well-known that a Hurwitz group of smallest order is a simple group, and as there's no Hurwitz group of order 84, then the discrete group of PSL(2,R) over the unique field with seven elements formed by the integers {0,1,2,3,4,5,6} under addition and multiplication modulo 7, has order twice 84 = 168, hence is "simple." This gives us the smallest finite projective space, the "Fano Plane" – which we navigate with the PSL(2,7) triangle! Given that all the "A, D, E" apparatus and form language can justly be said to be the extension of the symmetric wealth of implications stuffed inside the icosahedral E8, which in turn – as Galois showed – was a 60-element rotation group which is the first non-trivial group to be "simple," with the first Hurwitz group being the second, and Conway's "Monster" being the last, it is easy to believe two things: first, that the wealth of implications stuffed inside CDP ZD's, which are clearly as intimately connected to the Hurwitz group as "A, D, E" phenomena are to Galois' alternating group, is not just unexplored but huge; second, that both the Galois and Hurwitz "simple" groups are minuscule compared to Conway's, with so many others straddling the enormous gap between them, suggesting that all our mathematical understanding at present is likely next to nothing from the "God's eye" view! A good intro to Hurwitz groups, in a context which explicitly includes the general setting wherein Timaeus and hyperbolic (hence Dolgachev's, per the coming Note 20) triangles emerge, is Chapter 6 of Gareth A. Jones and David Singerman, <u>Complex Functions: An algebraic and geometric viewpoint</u> (Cambridge University Press, 1987), coming to a head for us on p. 261. (P.S.: the Pathions have 77 Box-Kites; add the Sedenion 7, and get 84 again!)

[16]  Algirdas Greimas, <u>Structural Semantics: An Attempt at a Method</u>, Daniele McDowell, Ronald Schleifer, Alan Velie, transl. (Lincoln NB and London: University of Nebraska Press, 1983; French original, 1966), p. 129. If Greimas's "two marks" be a pair of mutual ZD's, instead of oppositely signed instances of "X" (which we take to designate opposite termini of some Box-Kite edge); and if, finally, his "+" be replaced by a multiplicative sign, the correspondence between his setup and ours would be nigh complete. His choices of mathematical symbols being purely conventional (and hence, for the purposes of his argument, merely arbitrary), those suggested are just as viable as Greimas's own – and offer up

the bonus of framing a substantive theory, as we will soon see.  Note, too, we get the inherent split of his "marks" into dyads (the unit pairs comprising Assessors) "for free" this way – which gives us two mutually reinforcing readings of his "re-marks" simultaneously.

[17]    John Allen Paulos, Mathematics and Humor (Chicago and London:  University of Chicago Press, 1980) makes this point in great detail.  Its analysis, unfortunately, doesn't take things much further than the Cusp's immediate "big brother," the 3-control Swallow's Tail.

[18]    Jacques Derrida, Margins of Philosophy, Alan Bass, transl. (Chicago:  University of Chicago Press, 1982; French original, 1972), p.321.

[19]    This returns us to the theme of Note 15:  the Hurwitz simple group – which, after all, is seminal for us, and hence our paradigmatic instance of the Fuchsian groups which were the subject of Poincaré's "Aha!" experience.  Returning to Jones and Singerman's text cited in that note, we read on p. 221 that "The connection between the group PSL(2,R) and hyperbolic geometry was discovered by Henri Poincaré (1854-1912) and published in 1882.  (Poincaré used this discovery to illustrate the sometimes spontaneous nature of mathematical creativity.  He relates how, when boarding a bus, the idea suddenly came to him that the transformations he had used to define Fuchsian functions were identical with those of non-Euclidean geometry.)"

[20]    V. I Arnol'd, "Normal forms for functions near degenerate critical points, the Weyl groups of $A_k$, $D_k$, $E_k$ and Lagrangian singularities," Funct. Anal. Appl. 6 (1972), 254-272; transl. from Funkts. Anal. Prilozh 6:4 (1972), 3-25.  In a later paper translated in Russ. Math. Surv. 33:5 (1978), 99-116, Arnol'd extended the "A, D, E Correspondence" to "B, C, F" (on manifolds with boundary); later workers in his school have stirred up so many correspondences (to every item, in fact, in Coxeter's very long list) that it might be better to dub the ensemble "Alphabet Soup."  The first cross-disciplinary study of this "mysterious unity of all things," as it was described by the very unmystical Arnol'd, is M. Hazewinkel, W. Hesselink, D. Siersma, F. D. Veldkamp, "The ubiquity of Coxeter-Dynkin diagrams (an introduction to the A-D-E Problem)," Nieuw Arch. Wisk. 25 (1977), 257-307.  More recently, John Baez has devoted many installments of his weekly website offerings to this theme; see, too, the discussions and numerous links on Tony Smith's website.

[21]    See the discussion in V. I Arnol'd, S. M. Gusein-Zade and A. N. Varchenko, Singularities of Differentiable Maps, Vol. I:  The Classification of Critical Points, Caustics and Wave Fronts (Boston:  Birkhäuser, 1985), pp. 183-6, 242-6.  These 3 triangles are the bases of the 3 infinite families of "unimodular" germs, typically written with 3 subscripts to the letter "T" (one for each power its X, Y, or Z "behavior" displays in the germ), these subscripts being just the fraction of a straight angle included in successive angles of the triangle.  Hence, the Equilateral, Right Isosceles, and 30-60-90 of Plato's cosmogony are rendered in the following manner:  $T_{3,3,3}$; $T_{2,4,4}$; $T_{2,3,6}$ (also expressed, just to make things confusing, as $P_8$, $X_9$, and $J_{10}$ in the literature, with $X_9$ being Jim Callahan's preferred way of indicating the famous "Double Cusp" Catastrophe).  The "maximal strata" of these 3 are the finite E series of "simple" germs – and these, in turn, can be interpreted as the symmetries of the Platonic Solids:  $E_6$, then, represents the tetrahedron; $E_7$, the symmetrically dual cube and octahedron ; and $E_8$, the dodeca- and icosa- hedral duals.  These designations arise from a special "A-D-E" relationship known as the "McKay Correspondence," which is described in the just-cited reference (on p. 184) as follows:  "These singularities may also be obtained from the regular polyhedra in three-dimensional Euclidean space or more precisely from the discrete subgroups of the group SU(2):  they describe relations between the basic invariants of the groups.  $A_k$ corresponds to the polygons, $D_k$ to the dihedra (the two-sided polygons)" – and the E series, to the Solids as just indicated.  Arnol'd's associate Dolgachev has meanwhile (see p. 185) extended these amazingly simplistic yet deep correspondences into higher reaches:  the 14 "exceptional" forms among the unimodular germs (i.e., comprising finite listings, like the E series among the "elementary" germs), "may be obtained from the 14 triangles in the Lobechevskii plane [and hence realizable by M. C. Escher-style "hyperbolic tessellations"] or more precisely from the discrete subgroups of the group SU(1,1) determined by them" – with the triangles' angles, again, bearing exact correspondences to germinal expressions (in terms of what are now called the "Dolgachev numbers" linked to them).

[22]    Umberto Eco, A Theory of Semiotics (Bloomington IN:  Indiana University Press, 1976), p. 121.  Also cited in Schliefer, op. cit., p. 79.

[23]    Wolfram has a very useful discussion of this in the NKS pages on the "History of Chaos Theory," with Hadamard's 1898 thoughts, and Pierre Duhem's reflections on their significance in 1908, being referenced therein, on p. 971.

[24]    James Callahan, "Bifurcation Geometry of $E_6$," Mathematical Modeling, 1 (1980), 283-309 lays out the "Dynkinese" implicit in $E_6$ and how it decomposes into 3 classes of maximal strata:  Petitot's Dual

and Standard Butterfly forms, and the Umbilics which will concern us when we consider Trip-Sync-based "message-passing," in relation to the Semiotic Square as postmodernists see it. In sequal studies (the last two of which were co-authored with psychoanalyst Jerome Sashin), he developed the "tableau" approach as a method for modeling emotional affect response: see "A Geometric Model of Anorexia and Its Treatment," Behavioral Science, 27 (1982), 140-154; "Models of Affect-Response and Anorexia Nervosa," in S. H. Koslow, A. J. Mandell, M. F. Shlesinger, eds., Conference on Perspectives in Biological Dynamics and Theoretical Medicine, Annals N. Y. Acad. Sci., 504 (1986), 241-259; and, "Predictive Models in Psychoanalysis," Behavioral Science, 35 (1990), 241-259.

[25]   R. E. Cawagas, "Loops Embedded in Generalized Cayley Algebras of Dimension $2^t$, r ≥2," Int. J. Math. Math. Sci., 28:3 (2001), 181-187. To recycle Footnote 4 from my NKS 2004 paper, Prof. Cawagas' "Finitas" software, which he makes freely available, lets one calculate loops and other nonassociative structures rather easily in higher-dimensional contexts; he in fact found an equivalent manner of expressing my tabulations from BK1 in "loop language," which he presented in a 2004 paper to the National Research Council of the Philippines, "On the Structure and Zero Divisors of the Sedenion Algebra"; the correspondence which ensued after he communicated this to me proved delightfully engaging, and was the direct trigger for "BK3." On reading this, he further suggested I was engaging in "representation theory" – an intimation which surprised me at first, but then became an obsession … upon which much of this monograph can be blamed!

[26]   Dually, we could assume the reverse, adopting an "odometer-wheel reading" convention, in which case, per the earlier discussion of U- and L- index flows around a Box-Kite's hexagon, counterclockwise and clockwise change places in what follows … and the flow-reversal happens along the Tray-Rack edge *opposite* that specified. (But either way, obstacle-bypassing still involves detouring to the Zigzag Assessor perpendicular to the Tray-Rack in question.)

[27]   Arnol'd has talked about this "boundary singularity" in many places, but the best starting point is his book largely inspired by its amazing history as well as its surprising wealth of applications, Huygens and Barrow, Newton and Hooke , Eric J. F. Primrose, transl. (Basel, Boston, Berlin: Birkhäuser Verlag, 1990; Russian original, 1989). See especially Chapter 3: "From Evolvents to Quasicrystals," pp. 53-66, where largely overlooked mathematics of two centuries' past are dusted off, put in their own and current contexts, and found to be an endless fount of surprises and new insights. Bennequin's lecture to the 1984/5 Bourbaki seminar, "Caustiques mystiques," first revealed, as Arnol'd tells us in his introduction, "that the first textbook on analysis, written by l'Hôpital from the lectures of Johann Bernoulli, contains a representation of the manifold of irregular orbits of the Coxeter group $H_3$ (generated by reflexions in the planes of symmetry of an icosahedron). This representation appears there not in connection with the group of symmetries of the icosahedron, but as a result of investigations of evolvents of plane curves with a point of inflexion, investigations very close to those of Huygens (and possibly even carried out by him, although the first publication was apparently due to l'Hôpital). Illustrations appearing in recent works on the connection between the icosahedron and singularities of evolutes and evolvents and, it should be said, obtained by modern mathematicians not without difficulty and even with the help of computers, were already known at that time." What will intrigue us later is this: the icosahedral "group of symmetries" is based on the 15 median-flips plus identity, hence resides in 16-D, with the singularity appearing on its 4-D boundary. But the Sedenions likewise live in 16-D, and – as Moreno points out – have a 4-D "boundary" of ZD-free *lebensraum* to complement the 12-D Box-Kites need. More, consider the "words" and their readings that Quincunx lanyards can generate. As Sails can be either Zigzags or Trefoils, here we have Feet or Hands: If the "detour" is to the Zigzag Assessor, we get 10 possible starting points to read or write "/////\\\\\"; if we avoid the traffic-flow obstacle by the Vent Assessor instead, we get 10 strings from this: "/\\/\\/\\." Since scalar value along a diagonal doesn't matter, we are free to treat signing as binary "flow-orientation." And each Tray-Rack has its own Quincunx string-set. So 2 types x 3 Tray-Racks x 10 string-readings each x 2 flow-reversals = 120. Note further than Feet and Hands each break into two mirror-opposite 5-diagonal substrings, and that the 16-D singularity has (among other "preparations") a dynamic defined by a "double Swallow's Tail" (a $Z^5$ germ, with 3 controls on the 3 low monomials, Z running down the left or right slope of Pascal's Triangle, depending upon whether it be set to X or Y), with the 4-D "content" – the boundary singularity – unfolding the symmetric terms $1, xy, x^2y^2, x^3y^3$ in a so-called "open" Swallow's Tail. Finally, consider the workings of RNAi, which seems to be the master mechanism for turning genes off: it hunts for chunks of DNA made of mirror-image string-pairs; the Swallowtail, meanwhile, models self-destructive processes (Thom's "suicide catastrophe"). In Thom's 1971 model, which will come up more than once here, the 5-fold-symmetric "pentagram" unfolding effectively ends the "transmission" between Double

Cusps if we adapt it to the heuristic in BK1, which is sketchily alluded to in Section 14. All of which commentary must suffer in the limbo of these notes until some later date, when something more formally elegant (or simulation-based in its substantiveness) can be said about it.

28     Jean Petitot-Cocorda, Physique du Sens: De la théorie des singularités aux structures sémio-narratives (Paris: CNRS Editions, 1992), henceforth PS, pp. 388-394, discusses this in some detail, with generous graphical accompaniments, in the follow-up discussions of open problems upon completion of his "Semiotic Square" arguments. In particular, he reverts to a second "Greimas Glitch," framed in a 1973 paper called "Un problème de sémiotique narrative: les object de valeur." Here, Greimas focused on the conundrum posed by the seemingly simple (indeed, ubiquitous) phenomenon of "double transfer," which I must somehow magically negotiate, in spite of its epistemological difficulty, each time I un-pocket two quarters to procure the morning paper. This problem is so mysteriously intractible that "even if one makes use of the Cuspoid $X^\circledast$ [the Star], it is not possible to schematize the double transfer in a simple manner." As we'll see later, the way the Star relates to our "twist product" is via a dynamic interaction between two (or even three ) Double Cusps – something we "get for free" by "twisting," since the pairs of ZD's before and after it reside, as said above in the main text, in different Box-Kites.

29     Greimas's own "actantial" analysis clearly requires this double-articulated approach to his SS: on the one hand, Thom's CT models explicitly derive from the "actantial" approach of Tesnère's Éléments de syntaxe structurale (Paris: Klincksieck, 1959), also a key source of Greimas's own thoughts, as Petitot has discussed at some length. The Double Cusp's "generic sentence" with four minima has been an oft-discussed feature of any language models of the case-based variety since such models originated: see Zeeman, op. cit., where we are told in 1974 that the "most interesting key to linguistics, therefore, seems to lie in the study of paths in the double cusp, and the associated sequences of entrances, exits and transfer between the 4 actors involved, and the comparison of these paths with Thom's original classification of basic sentences." (p. 632) These remarks should be compared with Greimas's own, from his 1976 study of Maupassant, cited in Schleifer, op. cit., p. 156: "analytical experience – both our own and that of other semioticians – has convincingly demonstrated that, to account for texts even a little complex, it is necessary to consider the possibility of exploding (éclatement) any actant into at least four actantial positions" on the Semiotic Square.

30     Algirdas Greimas, op. cit., p. 19.

31     Ibid., p. 8.

32     In what follows, I depend directly upon Greimas and François Rastier, "The Interaction of Semiotic Constraints," Yale French Studies, 41 (1968), 86-105 – the essay anouncing the Semiotic Square and its isomorphism with (among other structures) the Klein Group (p. 88), and Lévi-Strauss's "canonical law of myths" (p. 89). My ZD labeling scheme serves to answer a deep question Petitot asks after showing the inadequacy of Boolean logic to handle the task: "How to conceive of relations of junction (conjunction/disjunction) and of reciprocal presupposition which are 'pure', that is to say independent of the specific semes whose value they determine?" (MM, p. 30)

33     Algirdas Julien Greimas and Jacques Fontanille, The Semiotics of Passions, Paul Perron and Frank Collins, transl. (Minneapolis and London: U. of Minnesota Press, 1993; French original, 1991). As the foreword to the English translation of this last book-length work by Greimas puts it, "The very fact that one of the first definitions in Descartes's theory deals with the relation between active nad passive explanations of passions had as a consequence the linking of passions not only to reason but also to actions. Within this framework Greimas's strategy was to work out a semiotic of passions, a pathemic semiotics founded on his own previously elaborated semiotics of action." (P. ix) The work concludes with a study of "the general syntax of jealousy," wherein Racine's Roxane "gives us an almost complete and particularly detailed realization of the microsequence" of abstract patterns of emotional investment and hedging that constitute the "fiduciary" dynamics of jealous entanglement as the Semiotic Square would model it (with the key characters at 3 of the Square's 4 corners, as well as giving us "a fourth position," encountered in the form of the 'gift of trust' that seems to correspond to an attachment formed by 'abandoning distrust.' (P.182, where the diagram appears)

34     The Chingon graphic for S=15 bears an obvious relationship to the combined effect of two Cesaro "triangle sweeps" initiated by [0,1] and [1,0], depicted on pp. 64-5 of Benoit Mandelbrot, The Fractal Geometry of Nature (San Francisco: W. H. Freeman, 1983; first edition, 1972). As we'll argue more formally later, the key is in the "spectrographic" effect of all bits to the left of that in the '4's' position: each row and column in the S=15 emanation table will have all entries for 8, 16, 24 filled (save for those which fall on the long diagonals); redoubling for the same S and keeping the screen size constant would thereby induce a doubling of the number of such lines (as well as new subdiagonals in quadrants, etc.).